\pdfoutput=1 
\documentclass[11pt,oneside,reqno]{amsart}

\usepackage[OT2,T1]{fontenc}
\usepackage[utf8]{inputenc}
\usepackage{amssymb,bm,bbm,mathtools}

\usepackage[dvipsnames]{xcolor}
\usepackage[
    colorlinks=true,
    linkcolor=Maroon,
    citecolor=JungleGreen,
    urlcolor=NavyBlue]{hyperref}

\usepackage{xcolor,colortbl}
\usepackage[capitalize,nameinlink]{cleveref}

\usepackage{enumitem}
\usepackage{dsfont}
\usepackage{tikz}
\usetikzlibrary{patterns}
\usepackage{subcaption}
\usetikzlibrary{arrows}
\usepackage{multirow}

\usepackage{lipsum}
\usepackage{soul}

\newtheorem{theorem}{Theorem}[section]

\newtheorem{lemma}[theorem]{Lemma}
\newtheorem{proposition}[theorem]{Proposition}
\newtheorem{corollary}[theorem]{Corollary}

\newtheorem{remark}[theorem]{Remark}

\numberwithin{equation}{section}

\DeclareMathOperator{\re}{Re}

\DeclareMathOperator{\im}{Im}

\DeclareMathOperator{\SO}{SO}
\DeclareMathOperator{\Spin}{Spin}
\DeclareMathOperator{\End}{End}
\DeclareMathOperator{\Id}{Id}
\DeclareMathOperator{\tr}{tr}

\DeclareMathOperator{\supp}{supp}

\DeclareMathOperator{\const}{const.}

\newcommand{\df}[1]{\textrm{d}#1}

\makeatletter
\def\@tocline#1#2#3#4#5#6#7{\relax
  \ifnum #1>\c@tocdepth 
  \else
    \par \addpenalty\@secpenalty\addvspace{#2}%
    \begingroup \hyphenpenalty\@M
    \@ifempty{#4}{%
      \@tempdima\csname r@tocindent\number#1\endcsname\relax
    }{%
      \@tempdima#4\relax
    }%
    \parindent\z@ \leftskip#3\relax
    \advance\leftskip\@tempdima\relax
    \rightskip\@pnumwidth plus4em \parfillskip-\@pnumwidth
    #5\leavevmode\hskip-\@tempdima
      \ifcase #1
       \or\or \hskip 2em \or \hskip 2em \else \hskip 3em \fi%
      #6\nobreak\relax
    \dotfill\hbox to\@pnumwidth{\@tocpagenum{#7}}\par
    \nobreak
    \endgroup
  \fi}
\makeatother

\usepackage[
      backend=biber,
      style=alphabetic,
      giveninits=true,
      maxbibnames=10,
      sorting=nyt,
      url=false,
      eprint=false,
      isbn=false
      ]{biblatex}

\DeclareFieldFormat[article,incollection]{title}{#1}
\DeclareFieldFormat{date}{#1}
\DeclareFieldFormat{pages}{{#1}}
\DeclareFieldFormat{doi}{
\href{https://www.doi.org/#1}{Doi}}
\DeclareFieldFormat{url}{
\href{#1}{Doi}}

\renewbibmacro{in:}{%
  \ifentrytype{article}{}{\printtext{\bibstring{in}\intitlepunct}}}

\begin{document}
\title[Dirac equation on spinor bundles]{Massless Dirac equation on spinor bundles  over real hyperbolic spaces}

\author{Long Meng, Hong-Wei Zhang, and Junyong Zhang}

\begin{abstract}
  We prove a sharp-in-time dispersive estimate of the Dirac equation on spinor bundles over the real hyperbolic space. Compared with the Euclidean counterparts, our result shows that the dispersive estimate differs between short and long times, reflecting the intuitive influence of negative curvature on the dispersion. Moreover, the well-known equivalence between dispersive estimates for Dirac and wave propagators in the Euclidean setting no longer holds in this context. This finding suggests that spinor fields are affected by the geometry at infinity of the manifold. As a key application, we establish an improved global-in-time Strichartz estimate, in the sense that there is no loss of angular derivatives and the admissible set is larger than previously known results in other settings.
\end{abstract}


\makeatletter
\@namedef{subjclassname@2020}{\textnormal{2020}
    \it{Mathematics Subject Classification}}
\makeatother

\maketitle
\section{Introduction}
This paper investigates the Dirac equation on the real hyperbolic space $\mathbb{H}^n$, which serves as a fundamental example of a negatively curved manifold. We consider the following time-dependent Cauchy problem:
\begin{equation}\label{eq: Dirac Intro}
    \begin{cases}
        (i\partial_t+D)U(t,x)=L(t,x),\\
        U(0,x)=U_0(x),
    \end{cases}
\end{equation}
where the massless Dirac operator $D$ is known to act on the smooth sections of the homogeneous vector bundles over $\mathbb{H}^n$. See the definitions \eqref{eq:D} in local coordinates or \eqref{eq:D2} through the Clifford multiplication and the covariant derivative.

The study of dispersive equations in non-flat settings have attracted considerable attention in recent years, particularly in the context of the Schrödin-ger, wave, and Klein-Gordon equations. Concerning the compact manifolds, we refer, for instance, to \cite{Bou93,BGT04,BHS20,Zha20}. Some work on hyperbolic spaces can be found in \cite{Tat01,Ban07,BCS08,AP09,IS09,MT11,AP14}, while references \cite{Zha21,AZ24} explore their higher rank generalizations to noncompact symmetric spaces. For related studies on asymptotically hyperbolic manifolds, spherically symmetric manifolds, and other negatively curved manifolds, see \cite{BD07,BGH10,Bou11,Che18,SSW19,SSWZ20,SSWZ22}. More details and references can be found therein.

Compared to the Schr\"odinger, wave, and Klein-Gordon equations, the study of the Dirac equation is more challenging, even in the Euclidean setting. See, for instance, \cite{EV97,BH15,BH16,DFS07}. In addition to the difficulty arising from the spin structure, the involvement of the covariant derivative makes the Dirac operator on non-flat manifolds more like the one with electromagnetic potentials on flat manifolds, which are hard to be analyzed since these are variable coefficient operators.

Nevertheless, some progress has been made recently in the study of the Strichartz estimates for the Dirac equation on non-flat manifolds. In \cite{CS22,BCSZ22,CSM24}, the authors studied the Strichartz estimates on asymptotically flat manifolds and certain spherically symmetric manifolds. Recently, using an improved WKB approximation, the Strichartz estimates with wave and Klein--Gordon admissible sets are also proved for the Dirac equation on compact manifolds without boundary in \cite{CDM23}. Compared with the Strichartz estimates in the Euclidean setting, there is a loss of regularity in the aforementioned results. Some of these losses are essential due to the compactness, while others arise from technical issues, such as the angular regularity loss discussed in \cite{CS22,BCSZ22}.

In this paper, we study the Dirac equation on the $n$-dimensional real hyperbolic space $\mathbb{H}^n$, which is a Riemannian manifold with constant negative sectional curvature $-1$, and also can be regarded as a homogeneous space $G/K$ with $G=\Spin_{e}(n,1)$ and $K=\Spin(n)$. In this setting, the Dirac operator is known to act on smooth spinors of the spinor bundles $\Sigma(G/K)$ over $\mathbb{H}^n$. Let $L^p(G,\tau)$ be the $L^p$ space of spinors and $\mathcal{R}=-n(n-1)$ be the normalized scalar curvature of $\mathbb{H}^n$. See Section \ref{sec: prelim} for more details about these notations. 

The Dirac operator has both positive and negative spectra. The wave propagates in different directions for the positive and negative spectral parts. To state the dispersive property, we introduce the absolute value of the Dirac operator, $|D|=\sqrt{D^* D}=\sqrt{D^2}$, since $D$ is an essentially self-adjoint operator on $\mathbb{H}^n$. Note that this does not result in any loss of generality (see Section \ref{sec: kernel}). Our first result is the following sharp-in-time dispersive estimate of the smoothed Dirac propagator associated with the operator $|D|$.

\begin{theorem}\label{main thm dispersive}
  Let $2<q\le\infty$ and $\theta\ge(n+1)(1/2-1/q)$. Then there exists a constant $C_{n,q}>0$ such that, for all $t\in\mathbb{R}^*$,
  \begin{align*}
    \left\|\left(D^2-\frac{\mathcal{R}}{4}
    \right)^{-\frac{\theta}{2}} e^{it|D|}
    \right\|_{L^{q'}(G,\tau)\rightarrow L^{q}(G,\tau)}\,
    \le\,C_{n,q}
    \begin{cases}
      |t|^{-(n-1)(\frac{1}{2}-\frac{1}{q})}
      &\,\,\,\textnormal{if}\quad |t|<1,\\
      |t|^{-1}
      &\,\,\,\textnormal{if}\quad |t|\ge1.
    \end{cases}
  \end{align*}
\end{theorem}

\begin{remark}
The dispersive property behaves differently for short and long times, which differs from the Euclidean results. This is consistent with the intuitive expectation that the dispersion is influenced by the underlying geometry, a phenomenon already observed for other dispersive equations in certain negatively curved settings. 
\end{remark}

\begin{remark}\label{rmk: dispersive}
For $|t|\leq 1$, the dispersive estimate is the same as for the wave equation in $\mathbb{R}^n$ or $\mathbb{H}^n$. However, the long-time dispersion of the Dirac propagator, which is sharp-in-time (see Remark \ref{rem:J0}), decays more slowly than that of the wave propagator on hyperbolic spaces, where the decay rate is $|t|^{-3/2}$, see \cite{AP14}. More precisely, unlike the scalar Laplace-Beltrami operator, which is independent of the spin structure, the spin angular momentum parameter $s=\pm 1/2$, combined with the covariant derivative, shifts the generalized eigenfunctions of the Dirac operator and alters the behavior of the Plancherel density $\mu(r)$ near the origin (see \eqref{mu}). Nevertheless, Theorem \ref{main thm dispersive} can be viewed as providing a stronger dispersive property compared to the Euclidean results, as the long-time decay rate is independent of the norm index $q$.
\end{remark}

This different dispersive estimate yields the following improved global-in-time Strichartz inequality on spinor bundles. For $s\in\mathbb{R}$ and $1<q<+\infty$, let $H^{s,q}(G,\tau)$ be the fractional Sobolev space defined by the Dirac operator, see Section \ref{sec: strichartz} for more details.

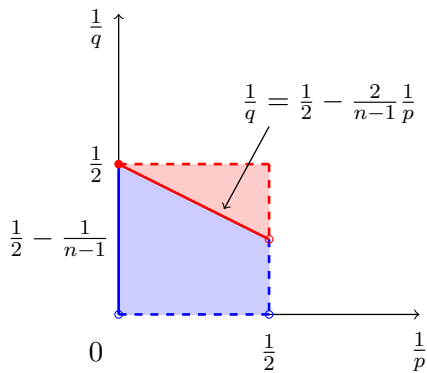
\begin{figure}
  \centering
  \begin{tikzpicture}
    \draw [->, line width=0.5] (2,0) -- (4,0);
    \draw [->, line width=0.5] (0,0) -- (0,4);

    \draw [blue, line width=0, opacity=0.2, fill=blue] (0,0) -- (0,2) -- (2,1) -- (2,0) -- cycle;
    \draw [blue, line width=1] (0,0) -- (0,2);
    \draw [blue, line width=1, dashed] (0,0) -- (2,0);
    \draw [blue, line width=1, dashed] (2,1) -- (2,0);

    \draw [red, line width=0, opacity=0.2, fill=red] (0,2) -- (2,2) -- (2,1) -- cycle;

    \draw [red, line width=1, dashed] (0,2) -- (2,2);
    \draw [red, line width=1, dashed] (2,2) -- (2,1);
    \draw [red, line width=1] (2,1) -- (0,2);

    \draw [red] (2,1) circle (0.05);
    \draw [red,fill=red] (0,2) circle (0.05);
    \draw [blue] (2,0) circle (0.05);
    \draw [blue] (0,0) circle (0.05);

    \node at (-0.3, -0.5) {$0$};
    \node at (4, -0.5) {$\frac{1}{p}$};
    \node at (-0.3, 3.8) {$\frac{1}{q}$};
    
    \node at (2, -0.5) {$\frac{1}{2}$};
    \node at (-0.3, 2) {$\frac{1}{2}$};
    \node at (-0.8, 1) {$\frac{1}{2}-\frac{1}{n-1}$};

    \draw [->, line width=0.5] (2,2.5) -- (1.4,1.4);
    \node at (2.8,2.8) {$\frac{1}{q}=\frac{1}{2}-\frac{2}{n-1}\frac{1}{p}$};
  \end{tikzpicture}
  \caption{Admissible set for $n\ge4$}
  \label{fig: admissible}
\end{figure}

\begin{theorem}\label{main thm strichartz}
  Let $(p_1,q_1)$ and $(p_2,q_2)$ be two admissible pairs, i.e., the points $(1/p_{1},1/q_{1})$ and $(1/p_{2},1/q_{2})$ belong to the admissible set (see Figure \ref{fig: admissible})
  \begin{align*}
    \left[0,\frac{1}{2}\right)\,
    \times\,
    \left(0,\frac{1}{2}\right)\,
    \bigcup\,
    \left\lbrace
    \left(0,\frac{1}{2}\right)
    \right\rbrace.
  \end{align*}
  Denote by 
  \begin{align*}
    \theta_i(p_i, q_i)\,
    =\,
    \frac{n+1}{2} \left( \frac{1}{2} - \frac{1}{q_{i}} \right)\,
    +\,
    \max \left\{ 0, \frac{n-1}{2} 
    \left( \frac{1}{2} - \frac{1}{q_i} \right) 
      - \frac{1}{p_i} \right\},
  \end{align*}
  for $i=1,2$, and let $\theta_i\ge\theta_i(p_i, q_i)$. Then, there exists a constant $C>0$, independent of $U$, $U_0$, and $L$, such that the solution of the linear Dirac equation \eqref{eq: Dirac Intro} satisfies the inequality
  \begin{align*}
    \|U\|_{L^{p_{1}}(\mathbb{R};\,H^{-\theta_{1},q_{1}}(G,\tau))}\,
    \le\,C\,
    \left(\|U_0\|_{L^{2}(G,\tau)}\,+\,
    \|L\|_{L^{p_{2}'}(\mathbb{R};\,H^{\theta_{2},q_{2}'}(G,\tau))}
    \right).
  \end{align*}
\end{theorem}

\begin{remark}
The Strichartz estimate in Theorem \ref{main thm strichartz} is sharp in regularity. Unlike the results in \cite{CS22,BCSZ22}, it does not involve any loss of angular regularity. Note also that the settings considered in \cite{BCSZ22} do not include $\mathbb{H}^n$.
\end{remark}

\begin{remark}\label{rem: admissible 1}
  The admissible set in Theorem \ref{main thm strichartz} corresponds to the points within the square shown in Figure \ref{fig: admissible} (the dashed lines and open points are not included). The red triangle represents the admissible set for the Strichartz estimate established via the $TT^*$ argument, while the blue trapezoid contains the admissible pairs obtained using the  Sobolev embedding. In the present context, the admissible set is larger: except for the red segment, the points inside the red triangle are typically not available in the Euclidean setting. This is due to our dispersive property, which differs from the Euclidean results, and is also in accordance with expectations in hyperbolic spaces. Conversely, the admissible set here is slightly smaller than that for the wave equation on hyperbolic spaces (see \cite[p. 966 Figure 1]{AP14})---the right edge of the square is not included in our setting. This is due to the slower large-time dispersion of the Dirac operator compared to the wave propagator, see Remark \ref{rmk: dispersive}. See also Remark \ref{rem: admissible 2} for figures of admissible sets in lower dimensions.
\end{remark}

\begin{remark}
  Understanding the dispersion of the propagator associated with the massive Dirac operator is a natural and interesting problem. However, its spectrum and the corresponding generalized eigenfunctions are more complicated, making the analysis more challenging. 
\end{remark}

\noindent\textbf{Organization of the paper.}
All the notations introduced in this section are explained in Section \ref{sec: prelim}, where we review the harmonic analysis tools on spinor bundles and the Dirac operator on hyperbolic spaces. Using these tools, we establish the sharp-in-time kernel estimate for the smoothed Dirac propagator in Section \ref{sec: kernel}. In the final section, we prove the main theorems presented in the Introduction. The dispersive estimate relies on a vector-valued Kunze-Stein phenomenon, while the Strichartz inequality follows from the classical duality argument combined with the Sobolev embedding.

Throughout this paper, the notation $A \lesssim B$ between two positive expressions indicates that there exists a constant $C>0$ such that $A \leq C B$. The notation $A \asymp B$  means that $A \lesssim B$ and $B \lesssim A$.

\section{Preliminaries}\label{sec: prelim}
In this section, we recall the Dirac operator on hyperbolic space. To do so, we first briefly review the homogeneous vector bundle over hyperbolic space. We also summarize the Fourier analysis on spinor bundles. We refer to \cite{Wal73,Kna86} and also to \cite{Cam97,CP01} for the following known facts.

\subsection{Real hyperbolic space}
For $n\ge2$, we denote the $n$-dimensional real hyperbolic space by 
\begin{align*}
  \mathbb{H}^{n}
  =\lbrace{x\in\mathbb{R}^{n+1}\,|\,
  x_{1}^{2}+\cdots+x_{n}^{2}-x_{n+1}^{2}=-1,\,
  x_{n+1}>0}\rbrace,
\end{align*}
which can be viewed as a symmetric space $\SO_{e}(n,1)/\SO(n)$ of noncompact type and rank one, and also as $G/K$ where $G=\Spin_{e}(n,1)$ and $K=\Spin(n)$ are the standard two-fold coverings of $\SO_{e}(n,1)$ and $\SO(n)$. Endowed with a suitable $G$-invariant Riemannian metric, the real hyperbolic space $\mathbb{H}^{n}$ is a Riemannian manifold with constant sectional curvature $-1$ and scalar curvature $-n(n-1)$.

Denote by $\mathfrak{g}=\mathfrak{spin}(n,1)\simeq\mathfrak{so}(n,1)$ and $\mathfrak{k}=\mathfrak{spin}(n)\simeq\mathfrak{so}(n)$ the corresponding Lie algebras of $G$ and $K$. According to the Cartan decomposition at the level of Lie algebra, one can write $\mathfrak{g}=\mathfrak{k}\oplus\mathfrak{p}$, where $\mathfrak{p}$ is identified with the tangent space of $\mathbb{H}^{n}$ at the origin. Let $\mathfrak{a}$ be the maximal abelian subspace in $\mathfrak{p}$ and $\mathfrak{a}^{+}$ be the corresponding positive Weyl chamber. Recall that $\mathfrak{a}\simeq\mathbb{R}$ and  $\mathfrak{a}^{+}\simeq\mathbb{R}^{+}$ in the present setting. We recall the following two decompositions of the group $G$, which play an important role in our study:
\begin{align*}
  \begin{cases}
    G=K(\exp\mathfrak{a})N
    &\qquad\textnormal{(Iwasawa)},\\
    G=K(\exp\overline{\mathfrak{a}^{+}})K
    &\qquad\textnormal{(Cartan)},
  \end{cases}
\end{align*}
where $N$ is a suitable nilpotent Lie group. The Haar measure $\df{x}$ on $G$ in the Cartan decomposition writes
\begin{align}\label{CartanHaar}
  \int_{G}f(x)\,\df{x}\,
  =\,
  \int_{K}\,\int_{0}^{\infty}\,(2\sinh r)^{n-1}\,\int_{K}
  f(k_{1}(\exp{r})k_{2}) \,\df{k_{2}}\,\df{r}\,\df{k_{1}},
\end{align}
where $\df{k_{1}}$ and $\df{k_{2}}$ are normalized Haar measures on the compact group $K$.

\subsection{Spinor bundle and spinors}
Given a unitary finite-dimensional representation $\tau$ of $K$ on a complex Hilbert space $V_{\tau}$, the associated homogeneous vector bundle over $G/K$ is denoted by
\begin{align*}
  \mathit{\Sigma}(G/K)\,=\,G\times_{K}V_{\tau},
\end{align*}
and its sections can be viewed as functions of type $\tau$ on $G$, i.e., as functions $f:G\rightarrow V_{\tau}$ which satisfy $f(xk)=\tau(k)^{-1}f(x)$ for all $x\in G$ and $k\in K$. We denote by $\mathcal{C}^{\infty}(G,\tau)$ the space of smooth $\tau$-type functions on $G$. We shall replace the $\mathcal{C}^{\infty}$ prefix with, e.g., $L^p$ or $\mathcal{S}$, for other function spaces. The homogeneous vector bundle over $\mathbb{H}^{n}$ is also called the spinor bundle and its sections are the spinors. We can write the spinor bundle as
\begin{align*}
  \mathit{\Sigma} \mathbb{H}^{n}=G\times_{K}V_{\tau_{n}},
\end{align*}
where $\tau_{n}$ is the unitary complex spin representation of $K$ on $V_{\tau_{n}}\simeq\mathbb{C}^{2^{[n/2]}}$. The representation $\tau_{n}$ is irreducible if the dimension $n$ is odd, but splits into two irreducible half-spin representations $\tau_{n}=\tau_{n}^{+}\oplus\tau_{n}^{-}$ if $n$ is even. In the latter case, we have 
\begin{align*}
  \mathit{\Sigma} \mathbb{H}^{n}=\mathit{\Sigma}^{+}\mathbb{H}^{n}\oplus\mathit{\Sigma}^{-}\mathbb{H}^{n}
  \qquad\textnormal{with}\qquad
  \mathit{\Sigma}^{\pm}\mathbb{H}^{n}=G\times_{K}V_{\tau_{n}^{\pm}}.
\end{align*}

Let $M$ be the centralizer of $\exp{\mathfrak{a}}$ in $K$. It is known that the boundary of the real hyperbolic space $\mathbb{H}^{n}$ is given by $K/M\simeq\mathbb{S}^{n-1}$. Let $\sigma_{n-1}$ be the complex spin representation of $M$ on $V_{\sigma_{n-1}}\simeq\mathbb{C}^{2^{[(n-1)/2]}}$. The space of spinors over $K/M$ can be viewed as the homogeneous vector bundle
\begin{align*}
  \mathit{\Sigma} \mathbb{S}^{n-1}=K\times_{M}V_{\sigma_{n-1}}.
\end{align*}

We denote by $\widehat{M}$ the unitary dual of $M$, which contains all equivalence classes of irreducible unitary representations of $M$. For any given irreducible unitary representation $\tau$ in $K$, let $\widehat{M}(\tau)$ be the set of $\sigma\in\widehat{M}$ occurring in the restrictions of $\tau$ on $M$. In our setting, we know that $\tau|_{M}$ is multiplicity free for any $\tau\in\widehat{K}$ and 
\begin{align}\label{tauM}
  \begin{cases}
    \tau_{n}|_{M}=\sigma_{n-1}^{+}\oplus\sigma_{n-1}^{-}
    &\qquad\textnormal{if $n$ is odd},\\
    \tau_{n}^{\pm}|_{M}=\sigma_{n-1}
    &\qquad\textnormal{if $n$ is even}.
  \end{cases}
\end{align}
See \cite[Propositions 6.2.3-6.2.4 and Theorems 8.1.3-8.1.4]{GW98}, see also \cite[Lemma 3.1]{CP01}.

\subsection{Dirac operator on hyperbolic space}\label{sec:2.3.1} We first begin with a brief overview of the construction of the Dirac equation in a non-flat (or non-Lorentzian) setting; we shall refer to \cite{CH96} for further details (see also Section 5.6 in \cite{PT09} and Section 2 in \cite{CS19}). Here for simplicity, we summarize the argument in \cite{CDM23}.

For any $n\geq2$ let $(\mathcal{M},g)$ be a $n$-dimensional Riemannian manifold endowed with a spin structure; then, the massless Dirac operator on $\mathcal{M}$ can be written as
\begin{align}\label{eq:D}
    D=-i\gamma^ae^j_{\; a}\nabla_j
\end{align}
and $\gamma^j$, $j=1,\dots,n$ is a set of matrices that satisfy the Clifford relation
\begin{align}\label{eq:Clifford-relation}
\gamma^h\gamma^j+\gamma^j\gamma^h=2\delta^{hj},\qquad h,j=1,\dots n.
\end{align}
The covariant derivative $\nabla_j$ for spinors is defined by
\begin{align}\label{conneq}
   \nabla_j=\partial_j+B_j,\quad j=1,2,\dots, n,
\end{align}
where $B_j$ writes
\[
B_j = \frac{1}{8} [\gamma^a,\gamma^b] \omega_j^{ab},
\]
and $\omega_j^{ab}$, called the {\em spin connection}, is given by
\begin{equation}\label{spincon}
    \omega_j^{ab}=e^h_{a}\partial_j e^{hb}+e_h^{a}\Gamma_{jk}^he^{kb},
\end{equation}
with the Christoffel symbol (or affine connection)
\begin{equation}\label{ichtus}
\Gamma_{jk}^h:=\frac{1}{2}g^{hl}(\partial_j g_{lk}+\partial_kg_{jl}-\partial_l g_{jk}),
\end{equation}
and the matrix bundle $e^h_{a}$ which is called {\em $n$-bein}:
\begin{equation}\label{dre}
    g^{hj} = e^h_{a}\delta^{ab}e^j_{b},
\end{equation}
where $\delta$ is the Kronecker symbol, and it connects the ``spatial'' metrics to the Euclidean one.

Finally, we precise the definition of matrices $\gamma^j$ which can be defined recursively(see, e.g., \cite{CH96} for more details). Before going further, we add subscript $n$ to the $\gamma^j$ matrices to keep track of the dimensions:
\begin{itemize}
\item {\em Case $n=2$.} We set
\begin{equation*}
\gamma_2^1=\left(\begin{array}{cc}0 &
i \\-i & 0\end{array}\right),\quad \gamma_2^2=\left(\begin{array}{cc}0 &
1 \\1 & 0\end{array}\right).
\end{equation*}
\item {\em Case $n=3$.} We set 
\begin{equation*}
\gamma_3^1=\gamma_2^1,\quad \gamma_3^2=\gamma_2^2,\quad \gamma_3^3=(-i)\gamma^1_2\gamma^2_2=\left(\begin{array}{cc}1 & 0 \\0 & -1\end{array}\right).
\end{equation*}
\item {\em Case $n>3$ even.} We set
\begin{equation*}
\gamma_n^j=\left(\begin{array}{cc}0 & i\gamma_{n-1}^j \\-i \gamma_{n-1}^j & 0\end{array}\right),\quad j=1,\dots,n-1,
\end{equation*}
and
\begin{align*}
     \gamma_n^n=\left(\begin{array}{cc}0 & I_{2^{\frac{n-2}2}} \\I_{2^{\frac{n-2}2}} & 0\end{array}\right).
\end{align*}
\item {\em Case $n>3$ odd.} We set
\begin{equation*}
\gamma_n^j=\gamma_{n-1}^j,\quad j=1,\dots,n-1,
\end{equation*}
and
\begin{align*}
    \gamma_n^n=i^{\frac{n-1}2}\gamma^1_{n-1}\cdot\dots\cdot\gamma^{n-1}_{n-1}=i^{\frac{n-1}2}\left(\begin{array}{cc}I_{2^{\frac{n-3}2}} & 0 \\0 &- I_{2^{\frac{n-{3} }2}}\end{array}\right).
\end{align*}
\end{itemize}

The Dirac operator in hyperbolic space is defined as a $G$-invariant differential operator acting on the sections of the spinor bundles. Let us recall its definition from \cite[Section 4.1]{CP01}. We denote by
\begin{align*}
  \gamma_{xK}\,
  :=\,T^*_{xK}(G/ K)\to \End \mathit{\Sigma}_{xK}(G/ K)
  \qquad\forall\,xK\in G/K = \mathbb{H}^{n}
\end{align*}
the extension of the Clifford multiplication to fibers of the spinor bundle, and by $\nabla:C^\infty(\mathit{\Sigma}(G/ K))\rightarrow C^\infty(T^*(G/ K)\otimes \mathit{\Sigma}(G/ K))$ the covariant derivative acting on smooth sections of $\Sigma(G/ K)$. Then the Dirac operator acting on $C^\infty(\mathit{\Sigma}(G/ K))$ can be defined as the composition of the covariant derivative with the Clifford multiplication, which has the local expression:
\begin{align*}
  D F(xK)\,
  =\,\sum_{j=1}^n\gamma_{xK}(e_j^*(xK))\,\nabla_{e_j}F(xK),
\end{align*}
where $(e_j)_{j=1}^{n}$ is any orthonormal frame on $G/K$ and $(e_j^*)^{n}_{j=1}$ is its dual co-frame. 

According to the Cartan decomposition $G=(\exp \mathfrak{p})K$, the Dirac operator $D$ is determined by its expression at the origin $eK$, then can be expressed as follows: if $F$ is a smooth $\tau$-type function, where either $\tau=\tau_{n}$ in odd dimension, or $\tau=\tau_{n}^{\pm}$ in even dimension, then for any $x\in G$,
\begin{align}\label{eq:D2}
    D F(x)\,
    =\,\sum_{j=1}^n \gamma_{eK}(e_j^*(eK))\,F(x:e_j),
\end{align}
where $(e_j)^n_{j=1}$ is an orthonormal basis of $\mathfrak{p}$ with dual basis $(e_j^*)_{j=1}^n$ and $\gamma_{eK}(e_j^*(eK))$ are matrices that satisfy the Clifford relation \eqref{eq:Clifford-relation}. Here we have used the following Harish-Chandra notation to denote the left and right differentiation of $F$ at $x$ with respect to the elements $X_{1},\dots,X_{\ell},Y_{1},\dots,Y_{m}$ in $\mathfrak{g}$:
\begin{align}\label{HC notation}
  F(X_{1}\cdots X_{\ell}&:x:Y_{1}\cdots Y_{m})\,
  =\,
  \frac{\partial}{\partial r_{1}}\Big|_{0}\,\cdots\,
  \frac{\partial}{\partial r_{\ell}}\Big|_{0}\,
  \frac{\partial}{\partial s_{1}}\Big|_{0}\,\cdots\,
  \frac{\partial}{\partial s_{m}}\Big|_{0}\notag\\
  &F((\exp r_{1}X_{1})\cdots(\exp r_{\ell}X_{1})\,x\,
  (\exp s_{1}Y_{1})\cdots(\exp s_{m}Y_{m})).
\end{align}

\subsection{Spherical Fourier transform on spinor bundles}
\label{sec:sphericalfourier}
In this subsection, we summarize the spherical Fourier transform on spinor bundles over $\mathbb{H}^n$. For more details about this part, see \cite[Sections 5 and 6]{CP01}. Given a finite-dimensional representation $(\tau,V_{\tau})$ of $K$, a function $F:G\rightarrow\End V_{\tau}$ is called $\tau$-radial on $G$ if it satisfies
\begin{align}\label{eq:tau-invar}
  F(k_{1}x k_{2})\,
  =\,\tau(k_{2})^{-1}F(x)\tau(k_{1})^{-1},
\end{align}
for all $x\in G$ and $k_{1},k_{2}\in K$. This notion extends the concept of radial functions on hyperbolic spaces to the case of spinor bundles. Note that if $F$ is $\tau$-radial, then the function $x\mapsto F(x)v$ is $\tau$-type for all $v\in V_{\tau}$. Moreover, a $\tau$-radial function is determined by its restriction to $\exp\overline{\mathfrak{a}^{+}}$ according to the Cartan decomposition. Since $\tau|_{M}$ is an irreducible unitary representation of multiplicity free, and $F|_{\exp\overline{\mathfrak{a}^{+}}}\in\End_{M}V_{\tau}$ where $M$ is the centralizer of $\exp{\mathfrak{a}}$ in $K$, Schur's lemma implies that $F|_{\exp\overline{\mathfrak{a}^{+}}}$ is scalar multiple of the identity on each $M$-submodule $V_{\sigma}$ of $V_{\tau}$. Then for each $\sigma\in\widehat{M}(\tau)$, there exists a scalar component $f_{\sigma}:\mathbb{R}_{+}\rightarrow\mathbb{C}$ of $F$ such that 
\begin{align}\label{eq:tau-invar2}
  F(\exp r)|_{V_{\sigma}}\,
  =\,
  f_{\sigma}(r) \Id_{V_{\sigma}}
  \qquad\forall r\ge0.
\end{align}
Therefore, together with the decomposition \eqref{CartanHaar}, the $L^2$ inner product of $\tau$-radial functions can be reduced as 
\begin{align*}
  (F,\widetilde{F})\,
  &=\,
  \int_{G}\tr\lbrace{F(x) \widetilde{F}(x)^{*}}\rbrace\,\df{x}\\
  &=\,
  \sum_{\sigma\in\widehat{M}(\tau)}\dim V_{\sigma}\,
  \int_{0}^{\infty}f_{\sigma}(r)\,\overline{\widetilde{f}_{\sigma}(r)}\,
  (2\sinh r)^{n-1}\df{r}.
\end{align*}

We denote by $\mathcal{S}(G,\tau,\tau)$ the Schwartz space for $\tau$-radial functions on $G$, i.e., the space of smooth $\tau$-radial functions satisfying
\begin{align*}
  \sup_{r>0}\,
  (1+r)^{N}\,e^{\frac{n-1}{2}r}\,
  \|F(D_{1}:\exp r:D_{2})\|_{\End{V_{\tau}}}\,
  <\infty,
\end{align*}
for any left-invariant differential operators $D_{1}$, $D_{2}$ on $G$ and $N\in\mathbb{N}$, where $F(D_{1}:\exp r:D_{2})$ is defined similarly to \eqref{HC notation}. The spherical Fourier transform of $F\in\mathcal{S}(G,\tau,\tau)$ is the collection of functions
\begin{align*}
  r\,\mapsto\,
  \lbrace{\mathcal{H}_{\sigma}^{\tau}(F)(r)
  }\rbrace_{\sigma\in\widehat{M}(\tau)},
\end{align*}
where $\mathcal{H}_{\sigma}^{\tau}(F)$ denotes the partial spherical Fourier transform of $F$:
\begin{align}\label{transform}
  \mathcal{H}_{\sigma}^{\tau}(F)(r)\,
  =\,
  (\dim V_{\tau})^{-1}\,
  \int_{G}\,\tr\lbrace{F(x)\Phi_{\sigma}^{\tau}(r,x^{-1})}\rbrace\,\df{x}.
\end{align}

Here $\Phi_{\sigma}^{\tau}$ denotes the vector-valued $\tau$-spherical functions on $G$, which serves a similar role as the exponential factor in Euclidean Fourier analysis. They are eigenfunctions for the algebra of left-invariant differential operators, particularly for the Dirac operator. For simplicity, we denote by 
\begin{align*}
  \Phi_{\pm}=\Phi_{\sigma_{n-1}^{\pm}}^{\tau_{n}}
  \qquad\textnormal{and}\qquad
  \Phi^{\pm}=\Phi_{\sigma_{n-1}}^{\tau_{n}^{\pm}}
\end{align*}
in odd and even dimensions respectively. According to \cite[Theorem 5.1]{CP01}, we know that
\begin{align}\label{eq:eigen}
  \begin{cases}
    D\Phi_{\pm}(r,\cdot)v\,
    =\,\pm r\Phi_{\pm}(r,\cdot)v
    &\qquad\textnormal{if $n$ is odd},\\
    D^{2}\Phi^{\pm}(r,\cdot)v^{\pm}\,
    =\,r^{2}\Phi^{\pm}(r,\cdot)v^{\pm}
    &\qquad\textnormal{if $n$ is even},
  \end{cases}
\end{align}
for all $r\in\mathbb{R}$, $v\in V_{\tau_{n}}$ in odd dimension, and $v^{\pm}\in V_{\tau_{n}^{\pm}}$ in even dimension. Both spherical functions $\Phi_{\pm}$ and $\Phi^{\pm}$ are identities at the origin of $G$, and holomorphic in the variable $r$. Moreover, they satisfy the following integral expressions:
\begin{align}\label{Phiodd}
  \Phi_{\pm}(r,x)\,
  =\,2\int_{K}
  e^{-(ir+\frac{n-1}{2})H(xk)}
  \tau_{n}(k)\circ P_{\sigma_{n-1}^{\pm}}\circ\tau_{n}(\underbar{$k$}(xk)^{-1})
  \,\df{k},
\end{align}
and
\begin{align}\label{Phieven}
  \Phi^{\pm}(r,x)\,
  =\,\int_{K}
  e^{-(ir+\frac{n-1}{2})H(xk)}
  \tau_{n}^{\pm}(k\underbar{$k$}(xk)^{-1})\,\df{k},
\end{align}
where \underbar{$k$}$(xk)$ and $H(xk)$ are respectively the $K$-component and $\mathfrak{a}$-component of $xk\in G$ in the Iwasawa decomposition, and $P_{\sigma_{n-1}^{\pm}}$ is a suitable projection of $V_{\tau_{n}}$ onto $V_{\sigma_{n-1}^{\pm}}$, see \cite[(4.10)]{CP01}.

According to \eqref{tauM}, in the following we denote by either $\tau=\tau_{n}$ and $\sigma=\sigma_{n-1}^{\pm}$ in odd dimension, or $\tau=\tau_{n}^{\pm}$ and $\sigma=\sigma_{n-1}$ in even dimension. For each such $\tau$ and $\sigma$, we define the function $\Pi_{\sigma}^{\tau}:K\times G\rightarrow\End V_{\tau}$ by 
\begin{align*}
  \Pi_{\sigma}^{\tau}(k,x)\,=\,
  \begin{cases}
    \Pi_{\sigma_{n-1}^{\pm}}^{\tau_{n}}(k,x)\,
    =\,2\tau_{n}(k)\circ P_{\sigma_{n-1}^{\pm}}\circ
    \tau_{n}(\underbar{$k$}(xk)^{-1})
    &\quad\textnormal{if $n$ is odd},\\
    \Pi_{\sigma_{n-1}}^{\tau_{n}^{\pm}}(k,x)\,
    =\,\tau_{n}^{\pm}(k\underbar{$k$}(xk)^{-1})
    &\quad\textnormal{if $n$ is even}.
  \end{cases}
\end{align*}
Then we can rewrite \eqref{Phiodd} and \eqref{Phieven} as 
\begin{align}\label{Phi}
  \Phi_{\sigma}^{\tau}(r,x)\,
  =\,\int_{K}
  e^{-(ir+\frac{n-1}{2})H(xk)}\,
  \Pi_{\sigma}^{\tau}(k,x)\,\df{k}.
\end{align}
Since $(\tau,V_{\tau})$ is a finite-dimensional representation of $K$, it follows that, for all $x\in G$ and $k\in K$,
\begin{align}\label{eq:estimate-Pi}
  \|\Pi_{\sigma}^{\tau}(k,x)\|_{\End{V_{\tau}}}\,
  =\,\mathrm{O}(1)
  \qquad\forall\,\sigma\in\widehat{M}(\tau).
\end{align}
It follows particularly that, for all $r\in\mathbb{R}$ and $x\in G$,
\begin{align}\label{Phi estimate}
  \|\Phi_{\sigma}^{\tau}(r,x)\|_{\End{V_{\tau_{n}}}}\,
  \lesssim\,\varphi_{0}(x),
\end{align}
where $\varphi_{0}$ denotes the scalar elementary spherical function defined by
\begin{align}\label{phi0}
  \varphi_{0}(x)\,
  =\,
  \int_{K}
  e^{-\frac{n-1}{2}H(xk)}\,\df{k}.
\end{align}
It is known that $\varphi_{0}$ is a bi-$K$-invariant function on $G$, i.e., $\varphi_{0}(x)=\varphi_{0}(\exp r)$, where $r\ge0$ is the middle component of $x\in G$ in the Cartan decomposition. Moreover, it satisfies
\begin{align}\label{phi0 estimate}
  \varphi_{0}(\exp r)\,
  \asymp\,
  (1+r)e^{-\frac{n-1}{2}r}
  \qquad\forall r\ge0.
\end{align}

We distinguish the notation of the partial spherical Fourier transform \eqref{transform} in odd and even dimensions such as 
\begin{align*}
  \mathcal{H}_{\pm}=\mathcal{H}_{\sigma_{n-1}^{\pm}}^{\tau_{n}}
  \qquad\textnormal{and}\qquad
  \mathcal{H}^{\pm}=\mathcal{H}_{\sigma_{n-1}}^{\tau_{n}^{\pm}}.
\end{align*}
In odd dimension, the transform $(\mathcal{H}_{+},\mathcal{H}_{-})$ is an isomorphism between the spaces $\mathcal{S}(G,\tau_{n},\tau_{n})$ and $\lbrace{(h_{+},h_{-})\in\mathcal{S}(\mathbb{R})^{2}\,|\,h_{\pm}(-r)=h_{\mp}(r)}\rbrace$. And the transform $\mathcal{H}^{\pm}$ is an isomorphism between $\mathcal{S}(G,\tau_{n}^{\pm},\tau_{n}^{\pm})$ and $\mathcal{S}(\mathbb{R})_{\textrm{even}}$. According to \cite[Theorem 6.3]{CP01}, the corresponding inverse formulae are given by
\begin{align*}
  F(x)\,
  =\,
  \int_{0}^{\infty}
  (\mathcal{H}_{+}(F)(r)\Phi_{+}(r,x)\,
  +\,\mathcal{H}_{-}(F)(r)\Phi_{-}(r,x))\,
  \mu(r)\,\df{r}
\end{align*}
in odd dimension, and 
\begin{align*}
  F^{\pm}(x)\,
  =\,
  \int_{0}^{\infty}
  \mathcal{H}^{\pm}(F^{\pm})(r)\Phi^{\pm}(r,x)\,
  \mu(r)\,\df{r}
\end{align*}
in even dimension. Here, the Plancherel density $\mu(r)$ is defined by 
\begin{align}\label{mu}
  \mu(r)\,
  =\,
  \begin{cases}
    \frac{\Gamma(\frac{n-1}{2})^{2}}{4\pi\Gamma(n-1)^{2}}\,
    \prod_{j=1}^{(n-1)/2}(r^{2}+(j-\frac12)^{2})
    &\qquad\textnormal{if $n$ is odd},\\
    \frac{1}{2^{2n-3}\Gamma(\frac{n}{2})^{2}}\,
    r\coth{(\pi r)}\,\prod_{j=1}^{n/2-1}(r^{2}+j^{2})  
    &\qquad\textnormal{if $n$ is even}.
  \end{cases}
\end{align}
From this explicit formula, we know that the density $\mu(r)$ satisfies, for all $n\ge2$ and $r\ge0$,
\begin{align}\label{mu estim}
  \mu(r)\,\lesssim\,(1+r)^{n-1}.
\end{align}
Note that the density $\mu(r)$ behaves as a constant around the origin, this behavior is different from the density $|\mathbf{c}(r)|^{-2}$ that occurs in the classical spherical Fourier transform in the scalar case: recall that $|\mathbf{c}(r)|^{-2}\sim r^{2}$ around the origin. Finally, we deduce directly from the calculations that 
\begin{align}\label{muk}
  \Big(\frac{\partial}{\partial r}\Big)^{k}\mu(r)\,
  =\,\mathrm{O}((1+r)^{n-1-k}),
\end{align}
for all $r\ge1$, $k\ge0$, and $n\ge2$. This means that the density $\mu(r)$ is an inhomogeneous differential symbol of order $n-1$.

\section{Dirac equation on hyperbolic spaces}\label{sec: kernel}
Throughout this section, we always denote by either $\tau=\tau_{n}$ and $\sigma=\sigma_{n-1}^{\pm}$ in odd dimension, or $\tau=\tau_{n}^{\pm}$ and $\sigma=\sigma_{n-1}$ in even dimension. Consider the time-dependent Dirac equation on the spinor bundle $G\times_{K}V_{\tau}$:
\begin{align}\label{Dirac0}
  \begin{cases}
    (i\partial_{t}\,+\,D)\,U(t,x)\,=\,0,\\
    U(0,x)\,=\,U_{0}(x),
  \end{cases}
\end{align}
where $U:\mathbb{R}^{*}\times G\rightarrow V_{\tau}$. Compared with previous studies of dispersive equations on hyperbolic spaces, the Dirac operator acts on smooth sections of the homogeneous vector bundle over $\mathbb{H}^{n}$ rather than on scalar functions, which complicated the arguments. In other words, we are considering vector-valued $\tau$-type functions instead of complex-valued functions on $\mathbb{H}^{n}$.

In addition, compared with the scalar Laplace–Beltrami operator on hyperbolic spaces, the eigenvalues of the Dirac operator $D$ on spinor bundles is less clear if $n$ is even (see \eqref{eq:eigen}). Nevertheless, according to spectral representation, we know from \eqref{eq:eigen} that
\begin{align*}
  |D|\,\Phi(r,\cdot)\,
  =\,|r|\,\Phi(r,\cdot),
\end{align*}
where either $\Phi=\Phi_{\pm}$ in odd dimension or $\Phi=\Phi^{\pm}$ in even dimension, and $|D|=\sqrt{D^{2}}$. So instead of studying directly the Dirac equation \eqref{Dirac0}, we consider the following half-wave type problem:
\begin{align}\label{Dirac eq}
  \begin{cases}
    (i\partial_{t}\,\pm\,|D|)\,U^\pm(t,x)\,=\,0,\\
    U^\pm(0,x)\,=\,U_{0}^\pm(x).
  \end{cases}
\end{align}

Before going further, let $P^+:= \mathbbm{1}_{[0,\infty)}(D)$ and $P^-:=\mathbbm{1}_{(-\infty,0)}(D)$ be the projections to the non-negative and negative spectrum of $D$ respectively. Formally,
\begin{align}\label{def: proj operator}
    P^\pm=  \frac{1}{2}\pm \frac{D}{2|D|},
\end{align}
with abusing the notation $s/|s|=1$ if $s=0$. Let $U_0^\pm=P^\pm U_0$  with $U_0$ being the initial data of the equation \eqref{Dirac0}, and note that 
\begin{align*}
    DP^\pm\,
    =\,
    \pm|D|P^\pm.
\end{align*}
Then it is easy to see that
\begin{align*}
    U^{+}+U^{-}\,=\,e^{it|D|}U_0^+\,+\,e^{-it|D|}U_0^-\,
    =\,e^{itD}\left(U_0^++U_0^-\right)\,
    =\,e^{itD}U_0\,=\,U.
\end{align*}

For further simplicity, we only consider the Dirac propagator $e^{it|D|}$ with $t\in\mathbb{R}^*$, and study the following equation
\begin{align}\label{Dirac eq'}
  \begin{cases}
    (i\partial_{t}+|D|)\,U(t,x)\,=\,0,\\
    U(0,x)\,=\,U_{0}(x).
  \end{cases}
\end{align}

Now consider the analytic family of operators $\lbrace\mathbf{D}_{t}^{\theta}\rbrace_{\theta\in\mathbb{C}}$, where 
\begin{align*}
  \mathbf{D}_{t}^{\theta}\,
  =\,\left(D^{2}-\frac{\mathcal{R}}{4}\right)^{-\frac{\theta}{2}}
  \,e^{it|D|}
\end{align*}
is the smoothed Dirac propagator acting on $\mathcal{C}^{\infty}(G,\tau)$. Recall that $\mathcal{R}=-n(n-1)$ is the normalized scalar curvature of $\mathbb{H}^{n}$. For simplicity, we denote $r_{0}=-\mathcal{R}/4$ in the following. Then $r_0$ is a positive constant depending only on the dimension $n$. 

Note that $U(t,x)=e^{it|D|}U_{0}(x)$ solves the Dirac equation \eqref{Dirac eq'}. According to \eqref{eq:eigen}, for each $\tau$ and $\sigma$, the partial spherical Fourier transform of $\mathbf{D}_{t}^{\theta}$ is given by $(r^{2}+r_{0})^{-\theta/2}e^{it|r|}$. By using the corresponding inverse formula, we write the convolution kernel of $\mathbf{D}_{t}^{\theta}$ as 
\begin{align*}
  K_{t}^{\theta}(x)
  =\sum_{\sigma\in\widehat{M}(\tau)}
  \int_{0}^{\infty}\,(r^{2}+r_{0})^{-\frac{\theta}{2}}\,e^{itr}\,
  \Phi_{\sigma}^{\tau}(r,x)\,\mu(r)\,\df{r},
\end{align*}
where the spherical function $\Phi_{\sigma}^{\tau}(r,x)$ and the Plancherel density $\mu(r)$ are defined in \eqref{Phi} and \eqref{mu}. See also \cite[(8.7)]{CP01} for a similar formula of the heat kernel on spinor bundles.

According to \eqref{tauM}, the sum over $\widehat{M}(\tau)$ has only a finite number of elements, then we focus in this section on establishing the pointwise estimate for the inner integral of the kernel $K_{t}^{\theta}(x)$. Let us divide it into two parts by using the cut-off functions $\chi_{0}(r)$ and $\chi_{\infty}(r)=1-\chi_{0}(r)$, where $\chi_{0}$ is compactly supported in the interval $[0,1)$ and $\chi_{0}=1$ in $[0,1/2]$:
\begin{align*}
  I_{0}(t,x)\,=\,\int_{0}^{\infty}
  \chi_{0}(r)\,(r^{2}+r_{0})^{-\frac{\theta}{2}}\,
  e^{itr}\,\Phi_{\sigma}^{\tau}(r,x)\,&\mu(r)\,\df{r},
\end{align*}
and
\begin{align*}
  I_{\infty}(t,x)\,=\,\int_{0}^{\infty}
  \chi_{\infty}(r)\,(r^{2}+r_{0})^{-\frac{\theta}{2}}\,
  e^{itr}\,\Phi_{\sigma}^{\tau}(r,x)\,&\mu(r)\,\df{r}.
\end{align*}
Note that the integrals $I_{0}$ and $I_{\infty}$ also depend on the representations $\tau$ and $\sigma$. We omit these notations for simplicity since all final estimates are uniform in $\tau$ and $\sigma$. 

Let us briefly summarize the strategy of the proofs. The integral $I_{0}$ is relatively easy to study, since the cut-off function $\chi_{0}$ guarantees its boundedness. The study of $I_{\infty}$ is more complicated. We will focus on the case where $\re\theta=(n+1)/2$, as this does not result in any loss of generality. On the one hand, all results obtained for the critical case $\re\theta=(n+1)/2$ clearly extend to the simpler case $\re\theta>(n+1)/2$. On the other hand, we will derive the dispersive estimates by interpolating the analytic family of operators $\lbrace\mathbf{D}_{t}^{\theta}\rbrace_{\theta\in\mathbb{C}}$ within the vertical strip $0<\re\theta<(n+1)/2$. The analysis of $I_{\infty}$ requires a careful examination of the corresponding oscillatory integrals to achieve sharp estimates. Especially in the subcase where $|x|/|t|>1/2$, we need a delicate analysis based on the Harish-Chandra expansion of the Jacobi functions. 

\begin{proposition}\label{estimate I0}
  Let $\theta\in\mathbb{C}$ and $\sigma\in\widehat{M}(\tau)$. Then the following estimate holds for all $t\in\mathbb{R}^{*}$ and $x\in G$:
  \begin{align}\label{estI01}
    |I_{0}(t,x)|\,\lesssim\,\varphi_{0}(x).
  \end{align}
  Moreover, if $|t|\ge1$, there exists $C_{n,\theta}>0$ such that 
  \begin{align}\label{estI02}
    |I_{0}(t,x)|\,\le\,C_{n,\theta}\,
    |t|^{-1}\,(1+|x|)\,\varphi_{0}(x).
  \end{align}
\end{proposition}

\begin{proof}
  According to the integral expression \eqref{Phi} of the spherical function, we write
  \begin{align*}
    I_{0}(t,x)\,
    &=\,\int_{0}^{\infty}\chi_{0}(r)\,(r^{2}+r_{0})^{-\frac{\theta}{2}}\,
    e^{itr}\,\Phi_{\sigma}^{\tau}(r,x)\,\mu(r)\,\df{r}\\
    &=\,\int_{K}e^{-\frac{n-1}{2}H(xk)}\,J_{0}(t,x,k)\,
    \Pi_{\sigma}^{\tau}(k,x)\,\df{k},
  \end{align*}
  where
  \begin{align}\label{def: int J0}
    J_{0}(t,x,k)\,
    =\,
    \int_{0}^{\infty}\chi_{0}(r)\,(r^{2}+r_{0})^{-\frac{\theta}{2}}\,
    e^{ir(t-H(xk))}\,\mu(r)\,\df{r}.
  \end{align}
  Since $\chi_{0}$ is compactly supported around the origin, it simply follows from \eqref{mu estim} that $J_{0}(t,x,k)=\mathrm{O}(1)$, and then by \eqref{eq:estimate-Pi},
  \begin{align*}
    |I_{0}(t,x)|\,\lesssim\,
    \int_{K}e^{-\frac{n-1}{2}H(xk)}\,
    |\Pi_{\sigma}^{\tau}(k,x)|\,\df{k}\,
    \lesssim\,\varphi_{0}(x).
  \end{align*}

  We immediately get estimate \eqref{estI02} from \eqref{estI01} in the case where $|x|/|t|>1/2$. Suppose that $|t|\ge1$ and $|x|/|t|\le1/2$. Using the integration by parts we obtain
  \begin{align*}
    J_{0}(t,x,k)\,
    &=\,\frac{\chi_{0}(r)\,\mu(r)}{i(t-H(xk))}\,
    (r^{2}+r_{0})^{-\frac{\theta}{2}}\,e^{ir(t-H(xk))}\Big|_{r=0}^{r=\infty}\\
    &-\frac{1}{i(t-H(xk))}\,\int_{0}^{\infty}e^{ir(t-H(xk))}\,
    \frac{\partial}{\partial{r}}
    \Big(\chi_{0}(r)\,(r^{2}+r_{0})^{-\frac{\theta}{2}}\,\mu(r)\Big),
  \end{align*}
  where the first term is $\mathrm{O}(|t-H(xk)|^{-1})$.

  According to \eqref{mu}, the last integral is bounded in odd dimension, since it vanishes for $r>1$. Note that the contribution of the term where the derivative is applied to $(r^{2}+r_{0})^{-\theta/2}$ is $\mathrm{O}(|\im\theta|)$. In even dimension, the factor $r\coth(\pi r)$ and its derivative are both bounded for $r\le1$. In fact, the hyperbolic sine function satisfies 
  \begin{align*}
    \sinh r\,\asymp\,\frac{r}{1+r}\,e^{r}
    \qquad\forall\,r\in\mathbb{R},
  \end{align*}
  from where we get that $\pi r/\sinh(\pi r)=\mathrm{O}(1)$ for all $r\in\supp\chi_{0}$. Since the hyperbolic cosine function is bounded around the origin, we get
  \begin{align*}
    r\coth(\pi r)\,=\,\mathrm{O}(1)
    \qquad\forall\,r\in\supp\chi_{0}.
  \end{align*}
  Moreover, we have
  \begin{align*}
    \Big|\frac{\partial}{\partial{r}}\Big(r\coth(\pi r)\Big)\Big|\,
    &=\,
    \Big|\coth(\pi r)\,-\,\frac{\pi r}{(\sinh(\pi r))^{2}}\Big|\\
    &\le\,
    |2\cosh{\pi r}|\,
    \Big|\frac{\pi r}{\sinh(\pi r)}\Big|\,
    \Big|\frac{1}{2\pi r}\,-\,\frac{1}{2\sinh(\pi r)\cosh(\pi r)}\Big|\\
    &=\,
    |2\cosh{\pi r}|\,
    \Big|\frac{\pi r}{\sinh(\pi r)}\Big|\,
    \Big|\frac{1}{2\pi r}\,-\,\frac{1}{\sinh(2\pi r)}\Big|,
  \end{align*}
  which is also $\mathrm{O}(1)$ for all $r\in\supp\chi_{0}$, since $r\mapsto1/(2\pi r)-1/\sinh(2\pi r)$ is a continuous function over $\mathbb{R}$ vanishing at the origin. 

  Therefore, we obtain 
  \begin{align*}
    |J_{0}(t,x,k)|\lesssim\,
    (1+|\im\theta|)\,
    |t-H(xk)|^{-1}.
  \end{align*}
  Since $|H(xk)|\le|x|$ for all $k\in K$ and $x\in G$, we know that $|t-H(xk)|^{-1}\le2|t|^{-1}$ in the present case where $|x|\le|t|/2$. Then we conclude that 
  \begin{align*}
    |I_{0}(t,x)|\,
    &\lesssim\,
    \int_{K}e^{-\frac{n-1}{2}H(xk)}\,
    |\Pi_{\sigma}^{\tau}(k,x)|\,|J_{0}(t,x,k)|\\
    &\lesssim\,
    (1+|\im\theta|)\,
    |t|^{-1}(1+|x|)\,\varphi_{0}(x),
  \end{align*}
  in the long-time case where $|t|\ge1$.
\end{proof}

\begin{remark}\label{rem:J0}
  Note that \eqref{def: int J0} is simply an oscillatory integral with a phase that has no critical points and an amplitude that is compactly supported near $0$, behaving as a constant within this support. This implies that the principal contribution of $J_0$ cannot decay faster than $|t|^{-1}$. See, for instance, \cite[p.334, Proposition 3]{Ste93}. Thus, the large-time decay in Proposition \ref{estimate I0} is sharp. 
\end{remark}

Studying $I^{\infty}(t,x)$ is more complicated, especially in the case where $r$ and $|x|/|t|$ are both large. Let us divide it again into two parts. Let $\chi_{0}(|t|r)$ and $\chi_{\infty}(|t|r)=1-\chi_{0}(|t|r)$ be two cut-off functions such that $\chi_{0}(|t|r)$ vanishes if $r\ge|t|^{-1}$ and $\chi_{\infty}(|t|r)$ vanishes if $r\le(2|t|)^{-1}$. We write $I^{\infty}=I_{-}^{\infty}+I_{+}^{\infty}$ with 
\begin{align*}
  I_{-}^{\infty}(t,x)\,=
  \,\int_{0}^{\infty}\chi_{0}(|t|r)\,\chi_{\infty}(r)\,
  (r^{2}+r_{0})^{-\frac{\theta}{2}}\,
  e^{itr}\,\Phi_{\sigma}^{\tau}(r,x)\,\mu(r)\,\df{r},
\end{align*}
and 
\begin{align*}
  I_{+}^{\infty}(t,x)\,=
  \,\int_{0}^{\infty}\chi_{\infty}(|t|r)\,\chi_{\infty}(r)\,
  (r^{2}+r_{0})^{-\frac{\theta}{2}}\,
  e^{itr}\,\Phi_{\sigma}^{\tau}(r,x)\,\mu(r)\,\df{r}.
\end{align*}
We do not need to find the estimates for all $\theta\in\mathbb{C}$, it is sufficient to find sharp estimates for $\theta\in\mathbb{C}$ with $\re\theta=(n+1)/2$ and then apply the interpolation in the vertical strip $0\le\re\theta\le(n+1)/2$ when we establish the dispersive properties. The estimate of $I_{-}^{\infty}$ follows easily, since it is defined over the interval $[0,|t|^{-1}]$.

\begin{proposition}\label{estimate Iinf-}
  Let $\theta\in\mathbb{C}$ with $\re{\theta}=\frac{n+1}{2}$. Then the estimate 
  \begin{align*}
    |I_{-}^{\infty}(t,x)|\,
    \lesssim\,
    \begin{cases}
      |t|^{-\frac{n-1}{2}}\,\varphi_{0}(x)\,
      &\qquad\textnormal{if}\,\,\,0<|t|<2,\\ 
      0
      &\qquad\textnormal{if}\,\,\,|t|\ge2,
    \end{cases}
  \end{align*}
  holds for all $x\in G$.
\end{proposition}

\begin{proof}
  Note at first that $\chi_{0}(|t|r)\chi_{\infty}(r)$ vanishes if $|t|\ge2$. Otherwise, the integral $I_{-}^{\infty}$ is defined over the interval $[1/2,1/|t|]$ with $|t|<2$. Proceeding as for Proposition \ref{estimate I0}, we write
  \begin{align*}
    I_{-}^{\infty}(t,x)\,
    =\,\int_{K}e^{-\frac{n-1}{2}H(xk)}\,J_{-}^{\infty}(t,x,k)\,
    \Pi_{\sigma}^{\tau}(k,x)\,\df{k},
  \end{align*}
  with
  \begin{align*}
    J_{-}^{\infty}(t,x,k)\,
    =\,
    \int_{0}^{\infty}\chi_{0}(|t|r)\,\chi_{\infty}(r)\,
    (r^{2}+r_{0})^{-\frac{\theta}{2}}\,
    e^{ir(t-H(xk))}\,\mu(r)\,\df{r}.
  \end{align*}
For $|t|<2$ and $\theta\in\mathbb{C}$ with $\re\theta=\frac{n+1}{2}$, we deduce from \eqref{mu estim} that 
\begin{align*}
  |J_{-}^{\infty}(t,x,k)|\,
  \lesssim\,
  \int_{\frac{1}{2}}^{\frac{1}{|t|}}
  r^{-\frac{n+1}{2}+n-1}\,\df{r}\,
  \lesssim\,|t|^{-\frac{n-1}{2}}.
\end{align*}
The proposition is therefore proved.
\end{proof}

Next, we divide the study of $I_{+}^{\infty}$ into two parts, depending on whether $|x|/|t|$ is large or small. If $|x|/|t|\le1/2$, we can get the desired results by repeating the integration by parts, and this condition ensures the estimate $|t-H(xk)|^{-1}\le2|t|^{-1}$. If $|x|/|t|>1/2$, the factor $|t-H(xk)|^{-1}$ may be uncontrollably large, and we need a more delicate argument.

\begin{proposition}\label{estimate Iinf+-}
  Let $\theta\in\mathbb{C}$ with $\re{\theta}=\frac{n+1}{2}$. Then, for all $t\in\mathbb{R}^{*}$ and $x\in G$ such that $|x|/|t|\le1/2$, there exists $C_{n,\theta}>0$ such that 
  \begin{align*}
    |I_{+}^{\infty}(t,x)|\,
    \le\,C_{n,\theta}\,\varphi_{0}(x)\,
    \begin{cases}
      |t|^{-\frac{n-1}{2}}
      &\qquad\textnormal{if $|t|<1$},\\
      |t|^{-N}
      &\qquad\textnormal{if $|t|\ge1$},
    \end{cases}
  \end{align*}
  for any $N\ge0$.
\end{proposition}

\begin{proof}
Let us denote by 
  \begin{align*}
    J_{+}^{\infty}(t,x,k)\,
    =\,
    \int_{0}^{\infty}\chi_{\infty}(|t|r)\,\chi_{\infty}(r)\,
    (r^{2}+r_{0})^{-\frac{\theta}{2}}\,
    e^{ir(t-H(xk))}\,\mu(r)\,\df{r},
  \end{align*}
  where, according to the definition of the cut-off functions,
  \begin{align}\label{cut-off}
    \chi_{\infty}(|t|r)\,\chi_{\infty}(r)\,
    =\,
    \begin{cases}
      \chi_{\infty}(|t|r)
      &\qquad\textnormal{if $|t|\le\frac{1}{2}$},\\
      \chi_{\infty}(r)
      &\qquad\textnormal{if $|t|\ge2$}.
    \end{cases}
  \end{align}

  In the long-time case where $|t|\ge2$, we have 
  \begin{align*}
    J_{+}^{\infty}(t,x,k)\,
    &=\,\lim_{\varepsilon\rightarrow0}
    \int_{0}^{\infty}\chi_{\infty}(r)\,
    (r^{2}+r_{0})^{-\frac{\theta}{2}}\,
    e^{ir(t-H(xk))-\varepsilon r}\,\mu(r)\,\df{r}\\
    &=\,\const\,\lim_{\varepsilon\rightarrow0}(i(t-H(xk))-\varepsilon)^{-N}\\
    &\times\int_{0}^{\infty}e^{ir(t-H(xk))-\varepsilon r}\,
    \Big(\frac{\partial}{\partial r}\Big)^{N}
    (\chi_{\infty}(r)\,(r^{2}+r_{0})^{-\frac{\theta}{2}}\,\mu(r))\,\df{r}
  \end{align*}
after $N$ times integrations by parts. Therefore, we obtain the following estimate, which is uniform in $\varepsilon$:
\begin{align*}
  |J_{+}^{\infty}(t,x,k)|\,
  \lesssim\,
  |t|^{-N}\,\int_{0}^{\infty}
  \Big|\Big(\frac{\partial}{\partial r}\Big)^{N}
    (\chi_{\infty}(r)\,(r^{2}+r_{0})^{-\frac{\theta}{2}}\,\mu(r))\Big|\,\df{r}.
\end{align*}
If there is a derivative that applies to the cut-off function $\chi_{\infty}(r)$, then the corresponding integral is reduced to the interval $[1/2,1]$, so it is bounded for all $\theta\in\mathbb{C}$. Otherwise, all derivatives are applied to the factor $(r^{2}+r_{0})^{-\theta/2}\,\mu(r)$, which is a differential symbol of order $-\theta+n-1$ according to \eqref{muk}. Therefore, for all $\theta\in\mathbb{C}$ with $\re\theta=(n+1)/2$ and $N>(n-1)/2$,
\begin{align*}
  |J_{+}^{\infty}(t,x,k)|\,
  \lesssim\,
  (1+|\im\theta|^{N})\,|t|^{-N}.
\end{align*}

When $|t|\le1/2$, we have similarly
\begin{align*}
  |J_{+}^{\infty}(t,x,k)|\,
  \lesssim\,
  |t|^{-N}\,\int_{0}^{\infty}
  \Big|\Big(\frac{\partial}{\partial r}\Big)^{N}
    (\chi_{\infty}(|t|r)\,(r^{2}+r_{0})^{-\frac{\theta}{2}}\,\mu(r))\Big|\,\df{r}.
\end{align*}
Suppose there are $N_{1}$ derivatives that apply to the cut-off function $\chi_{\infty}(|t|r)$ and $N_{2}$ derivatives that apply to the remaining factor, where $N_{1}+N_{2}=N$. If $N_{1}=0$, then 
\begin{align*}
  |J_{+}^{\infty}(t,x,k)|\,
  &\lesssim\,
  (1+|\im\theta|^{N})\,|t|^{-N}\,
  \int_{|t|^{-1}}^{\infty}
  r^{-\re\theta+n-1-N}\,\df{r}\\
  &\lesssim\,
  (1+|\im\theta|^{N})\,|t|^{-\frac{n-1}{2}},
\end{align*}
provided that $\theta\in\mathbb{C}$ with $\re\theta=(n+1)/2$ and $N>(n-1)/2$. If $N_{1}\ge1$, the integral is reduced to the interval $[|t|^{-1},2|t|^{-1}]$, and we obtain 
\begin{align*}
  |J_{+}^{\infty}(t,x,k)|\,
  &\lesssim\,
  (1+|\im\theta|^{N_{2}})\,|t|^{-N}\,|t|^{N_{1}}\,
  \int_{|t|^{-1}}^{2|t|^{-1}}
  r^{-\re\theta+n-1-N_{2}}\,\df{r}\\
  &\lesssim\,
  (1+|\im\theta|^{N})\,|t|^{-\frac{n-1}{2}},
\end{align*}
for $\theta\in\mathbb{C}$ with $\re\theta=(n+1)/2$.

When $1/2<|t|<2$, we can do the same integrations by parts as in the long-time case to make the integral $J_{+}^{\infty}$ converge, and add the factor $|t|^{-N}$ for any $N\ge0$, since it is nothing but a constant in this case.

We conclude that, for any $\theta\in\mathbb{C}$ with $\re\theta=(n+1)/2$ and any $N\ge0$,
\begin{align*}
  |I_{+}^{\infty}(t,x)|\,
  &\le\,\int_{K}e^{-\frac{n-1}{2}H(xk)}\,|\Pi_{\sigma}^{\tau}(k,x)|\,
  |J_{+}^{\infty}(t,x,k)|\,\df{k}\\
  &\le\,C_{n,\theta}\,\varphi_{0}(x)\,
  \begin{cases}
    |t|^{-\frac{n-1}{2}}
    &\qquad\textnormal{if $|t|<1$},\\
    |t|^{-N}
    &\qquad\textnormal{if $|t|\ge1$},
  \end{cases}
\end{align*}
where $C_{n,\theta}=1+|\im\theta|^{N}$.
\end{proof}

It remains for us to study the integral
\begin{align*}
  I_{+}^{\infty}(t,x)\,=
  \,\int_{0}^{\infty}\chi_{\infty}(|t|r)\,\chi_{\infty}(r)\,
  (r^{2}+r_{0})^{-\frac{\theta}{2}}\,
  e^{itr}\,\Phi_{\sigma}^{\tau}(r,x)\,\mu(r)\,\df{r}
\end{align*}
in the case where $|x|/|t|>1/2$. Note that if we do a similar integration by parts as above, we will obtain an uncontrollably large factor $|t-H(xk)|^{-1}$. To get the desired sharp estimates, we need a more delicate analysis for the scalar components of the vector-valued spherical function. These scalar components are defined as Jacobi functions and satisfy the Harish-Chandra expansion. Let us recall this expansion in the next paragraph.

\vspace{10pt}\noindent{\bf Harish-Chandra expansion of Jacobi functions.}
We adopt the notation used in \cite[Section 2.1]{Koo84} and consider only the Jacobi function $\phi_{r}^{(\alpha,\beta)}$ with special parameters $\alpha=n/2-1$ and $\beta=n/2$. The Jacobi function $\phi_{r}^{(n/2-1,n/2)}$ is the unique smooth even solution on $\mathbb{R}$ of the differential equation 
\begin{align}\label{Jacobi eq}
  \Big(\frac{\partial^{2}}{\partial s^{2}}\,
  +\,\big((n-1)\coth{s}\,
  +\,(n+1)\tanh{s}\big)\,
  \frac{\partial}{\partial s}\,
  +\,r^{2}\,+n^{2}\Big)v(s)\,
  =\,0
\end{align}
such that $v(0)=1$. The Jacobi function can be expressed by using the Gauss hypergeometric function as
\begin{align*}
  \phi_{r}^{(\frac{n}{2}-1,\frac{n}{2})}(s)\,
  =\,
  {}_{2}F_{1}\Big(\frac{n-ir}{2},\frac{n+ir}{2};\frac{n}{2};-(\sinh{s})^{2}\Big),
\end{align*}
and satisfies the following expansion for $r\notin i\mathbb{Z}$ and $s>0$:
\begin{align}\label{HC expansion1}
  \phi_{r}^{(\frac{n}{2}-1,\frac{n}{2})}(s)\,
  =\,\mathbf{c}(r)\,\Phi_{r}(s)\,+\,\mathbf{c}(-r)\,\Phi_{-r}(s).
\end{align}
Here the Harish-Chandra $\mathbf{c}$-function is defined by 
\begin{align}\label{cfunction}
  \mathbf{c}(r)\,
  =\,
  \frac{2^{n-ir}\,\Gamma(\frac{n}{2})\,
  \Gamma(ir)}{
    \Gamma(\frac{n+ir}{2})\,
    \Gamma(\frac{ir}{2}),
  } 
\end{align}
and 
\begin{align*}
  \Phi_{r}(s)\,
  =\,(2\cosh{s})^{ir-n}\,
  {}_{2}F_{1}\Big(\frac{n-ir}{2},1-\frac{ir}{2};1-ir;(\cosh{s})^{-2}\Big),
\end{align*}
is also a solution to the equation \eqref{Jacobi eq}. See, for instance, \cite[(2.4)--(2.17)]{Koo84}. Writing formally 
\begin{align}\label{HC expansion2}
  \Phi_{r}(s)\,
  =\,e^{(ir-n)s}\,\sum_{m=0}^{\infty}\Gamma_{m}(r)\,e^{-ms},
\end{align}
and inserting it in the differential equation \eqref{Jacobi eq}, it follows that the coefficients $\Gamma_{m}$ satisfy $\Gamma_{0}(r)=1$, $\Gamma_{2m-1}=0$, and the recurrence formula
\begin{align}\label{Gamma2m}
  \Gamma_{2m}(r)\,
  &=\,
  \sum_{m'=0}^{m-1}
  (\delta_{m'}^{m}(n+1)-1)\,
  \frac{2m'+n-ir}{m(m-ir)},
\end{align}  
where $\delta_{m'}^{m}=1$ if $m'$ and $m$ have the same parity, or $\delta_{m'}^{m}=0$ otherwise. See \cite[p.151]{FJ72}. \footnote{Note that the parameters $p$ and $q$ in Flensted-Jensen's paper satisfy $p=2\alpha-2\beta$ and $q=2\beta+1$, where $\alpha$ and $\beta$ are the notations used in Koornwinder's article \cite{Koo84}. In both works the parameter $\rho=p/2+q=\alpha+\beta+1$. In the present case where $\alpha=n/2-1$ and $\beta=n/2$, their $\rho$ is just $n$.}
The classical results in \cite{HC58,ST78} showed that, for any $m\in\mathbb{N}$, the factor $\Gamma_{m}(r)$ grows at most polynomially in $m$. Here we need a more precise estimate, which has the decay in $r$ at the same time. To see this, it is sufficient to rewrite
\begin{align}\label{Gamma decomp}
  \Gamma_{2m}(r)\,
  =\,
  \underbrace{\sum_{m'=0}^{m-1}
  \frac{\delta_{m'}^{m}(n+1)-1}{m}}_{:=\,C_{n,m}}
  +\,
  \underbrace{\sum_{m'=0}^{m-1}
  \Big(\frac{\delta_{m'}^{m}(n+1)-1}{m}\Big)\,
  \Big(\frac{2m'+n-m}{m-ir}\Big)}_{:=\,\widetilde{\Gamma}_{2m}(r)}.
\end{align}  
Then there exist $\nu_{1},\nu_{2}>0$ such that 
\begin{align}\label{Gamma decomp estim}
  C_{n,m}\,=\,\mathrm{O}(m^{\nu_{1}})
  \qquad\textnormal{and}\qquad
  \widetilde{\Gamma}_{2m}(r)\,=\,\mathrm{O}(m^{\nu_{2}}r^{-1}),
\end{align}
for all $m\in\mathbb{N}^{*}$ and $r>0$.

On the other hand, we need information about the Harish-Chandra $\mathbf{c}$-function. By using Legendre's duplication formula 
\begin{align*}
  \Gamma(2z)\,=\,2^{2z-1}\,\pi^{-\frac12}\,\Gamma(z)\,\Gamma(z+\frac12),
\end{align*}
see e.g. \cite[p.5, (15)]{EMOT53}, we can simplify the $\mathbf{c}$-function as 
\begin{align*}
  \mathbf{c}(2r)\,
  =\,
  \frac{2\,\Gamma(n-1)}{\Gamma(\frac{n-1}{2})}\,
  \frac{\Gamma(ir+\frac{1}{2})}{\Gamma(ir+\frac{n}{2})}.
\end{align*}
Together with the Stirling formula
\begin{align*}
  \Gamma(z)\,=\,
  (2\pi)^{-\frac{1}{2}}\,z^{z-\frac{1}{2}}\,e^{-z}\,
  (1+\mathrm{O}(|z|^{-1})),
\end{align*}
see \cite[p.47, (2)]{EMOT53}, we know that there exists a constant $C_{n}>0$ such that 
\begin{align}\label{c expansion}
  \mathbf{c}(2r)^{-1}\,
  &=\,
  \frac{\Gamma(\frac{n-1}{2})}{2\,\Gamma(\frac{n}{2})}\,
  e^{-\frac{n-1}{2}}\,
  \Big(\frac{ir+\frac{n}{2}}{ir+\frac{1}{2}}\Big)^{ir}\,
  \Big(ir+\frac{n}{2}\Big)^{\frac{n-1}{2}}\,
  (1+\mathrm{O}(|r|^{-1}))\notag\\
  &=\,
  C_{n}\,r^{\frac{n-1}{2}}\,(1+\mathrm{O}(|r|^{-1})).
\end{align}
The Stirling formula tells us that $\mathbf{c}(2r)^{-1}$ is an inhomogeneous differential symbol of order $(n-1)/2$. More precisely, since $\mathbf{c}(2r)=\mathbf{c}_{0}(r-i/2)$, where $\mathbf{c}_{0}$ denotes the classical Harish-Chandra function in the scalar case, see for instance \cite[p. 296]{Ion00}. It follows from \cite[(A.2)]{Ion00} that 
\begin{align}\label{c estimate}
  \Big|\Big(\frac{\partial}{\partial r}\Big)^{k}(r^{-1}\mathbf{c}(2r)^{-1})\Big|\,
  \lesssim\,(1+|r|)^{\frac{n-3}{2}-k},
\end{align}
for all $k\in\mathbb{N}$ and $r\in\mathbb{R}$ such that $|r|\ge1$.

Now let us turn to the study of $I_{+}^{\infty}$ in the remaining case where $|x|/|t|>1/2$. We will first consider the case of even dimensions where the scalar components have simpler expressions. The case of odd dimensions will be discussed later. Recall that the condition $|x|/|t|>1/2$ allows us to add the factor $|t|^{-N}(1+|x|)^{N}$, for any $N\ge0$. The polynomial growth in $|x|$ is harmless, since we have an exponential decay at the same time.

\begin{proposition}\label{estimate Iinf+even}
  Let $n$ be even and $\theta\in\mathbb{C}$ with $\re{\theta}=(n+1)/2$. Then, for all $x\in G$ and $t\in\mathbb{R}^{*}$ such that $|x|/|t|>1/2$, there exists $C_{n,\theta}>0$ such that 
  \begin{align*}
    |I_{+}^{\infty}(t,x)|\,
    \le\,C_{n,\theta}\,e^{-\frac{n-1}{2}|x|}.
  \end{align*}
\end{proposition}

\begin{proof}
Recall that a $\tau$-radial function $F:G\rightarrow\End V_{\tau}$ is completely determined by its restriction to the middle component $\exp\overline{\mathfrak{a}^{+}}$ of $G$ in the Cartan decomposition. For each $\sigma\in\widehat{M}(\tau)$, there exists a scalar component $f_{\sigma}$ of $F$ such that 
\begin{align*}
  F(\exp s)|_{V_{\sigma}}\,
  =\,
  f_{\sigma}(s) \Id_{V_{\sigma}}
  \qquad\forall\,s\ge0.
\end{align*}
See Section \ref{sec:sphericalfourier} for more details. In \cite[Theorem 5.4]{CP01} it was proved that the scalar component $\varphi^{\pm}(r,\cdot)$ of $\Phi^{\pm}(r,\cdot)=\Phi_{\sigma_{n-1}}^{\tau_{n}^{\pm}}(r,\cdot)$ is given by 
\begin{align}\label{scalar even}
  \varphi^{\pm}(r,s)\,
  =\,\Big(\cosh\frac{s}{2}\Big)\,
  \phi_{2r}^{(\frac{n}{2}-1,\frac{n}{2})}\Big(\frac{s}{2}\Big)
\end{align}
in even dimension. Since
\begin{align*}
  I_{+}^{\infty}(t,x)\,=
  \,\int_{0}^{\infty}\chi_{\infty}(|t|r)\,\chi_{\infty}(r)\,
  (r^{2}+r_{0})^{-\frac{\theta}{2}}\,
  e^{itr}\,\Phi_{\sigma}^{\tau}(r,x)\,\mu(r)\,\df{r}
\end{align*}
is also a $\tau$-radial function on $G$, we deduce from \eqref{scalar even} and \eqref{HC expansion1} that the scalar component of $I_{+}^{\infty}(t,\cdot)$ is given by 
\begin{align}\label{scalar I}
  \Big(\cosh\frac{s}{2}\Big)
  \int_{0}^{\infty}&\chi_{\infty}(|t|r)\,\chi_{\infty}(r)\,
  (r^{2}+r_{0})^{-\frac{\theta}{2}}\,e^{itr}\notag\\
  &\times\,\Big(\mathbf{c}(2r)\,\Phi_{2r}\Big(\frac{s}{2}\Big)\,
  +\,\mathbf{c}(-2r)\,\Phi_{-2r}\Big(\frac{s}{2}\Big)
  \Big)\,\mu(r)\,\df{r}.
\end{align}
According to \cite[Lemma 6.4]{CP01}, we know 
\begin{align*}
  \mu(r)\,\df{r}\,=\,2\pi^{-1}\,|\mathbf{c}(2r)|^{-2}\,\df{r},
\end{align*}
for all $r\in\mathbb{R}$. Since $\Gamma(\bar{z})=\overline{\Gamma(z)}$ for all $z\in\mathbb{C}$, we have the formula $|\mathbf{c}(2r)|^{-2}=\mathbf{c}(2r)^{-1}\mathbf{c}(-2r)^{-1}$. Together with the fact that $\mu$ is an even function, we get
\begin{align}\label{densities}
  \mathbf{c}(\pm2r)\,\mu(r)\,\df{r}\,
  =\,2\pi^{-1}\,\mathbf{c}(\mp2r)^{-1}\,\df{r},
\end{align}
for all $r>0$. Therefore, we deduce from \eqref{scalar I},  \eqref{densities}, and \eqref{HC expansion2} that the following estimate holds for all $t\in\mathbb{R}^{*}$ and $x\in\ G$:
\begin{align*}
  |I_{+}^{\infty}(t,x)|\,
  \lesssim\,
  e^{-\frac{n-1}{2}|x|}\,\sum_{m=0}^{\infty}\,
  e^{-\frac{m|x|}{2}}\,|J_{m}(t,|x|)|,
\end{align*}
where $J_{m}$ is defined in the next lemma, and satisfies the estimate \eqref{J estimate} below. Then the proposition is proved.
\end{proof}

\begin{lemma}
  Let $m\in\mathbb{N}$ and $J:\mathbb{R}^{*}\times\mathbb{R}_{+}^{*}\rightarrow\mathbb{C}$ be an oscillatory integral defined by
  \begin{align*}
    J_{m}(t,s)\,
    =\,
    \int_{0}^{\infty}\chi_{\infty}(|t|r)\,\chi_{\infty}(r)\,
    a_{m}^{\pm}(r)\,e^{ir(t\pm s)}\,\df{r},
  \end{align*}
  with amplitude $a_{m}^{\pm}(r)=(r^{2}+r_{0})^{-\theta/2}\Gamma_{m}(\pm2r)\mathbf{c}(\mp2r)^{-1}$. Then for $t\in\mathbb{R}^{*}$ and $s\in\mathbb{R}_{+}^{*}$ satisfying $s/|t|>1/2$, there exist constants $\nu>0$ and $C_{n,\theta}>0$ such that 
  \begin{align}\label{J estimate}
    |J_{m}(t,s)|\,
    \lesssim\,C_{n,\theta}\,(1+m^{\nu}),
  \end{align}
  for all $m\in\mathbb{N}$, provide that $\theta\in\mathbb{C}$ with $\re{\theta}=\frac{n+1}{2}$.
\end{lemma}

\begin{proof}
According to \eqref{Gamma decomp}, for all $m\geq 1$, we divide the amplitude $a_{m}^{\pm}$ as 
\begin{align*}
  a_{m}^{\pm}(r)\,
  =\,
  \underbrace{C_{n,m}\,(r^{2}+r_{0})^{-\theta/2}\,\mathbf{c}(\mp2r)^{-1}
  }_{=\,C_{n,m}\,a_{0}^{\pm}(r)}\,
  +\,
  \underbrace{(r^{2}+r_{0})^{-\theta/2}\,
    \widetilde{\Gamma}_{2m}(\pm2r)\,\mathbf{c}(\mp2r)^{-1}
  }_{:=\,\widetilde{a}_{m}^{\pm}(r)}.
\end{align*}
We know from \eqref{Gamma decomp estim} that $C_{n,m}$ is $\mathrm{O}(m^{\nu_{1}})$, and $\widetilde{\Gamma}_{2m}(\pm2r)$ is $\mathrm{O}(m^{\nu_{2}}r^{-1})$ for some constants $\nu_{1},\nu_{2}>0$. Together with \eqref{c estimate}, we deduce that the factor $\widetilde{a}_{m}^{\pm}(r)$ is a differential symbol of order $-2$, provided that $\re\theta=(n+1)/2$. Then the corresponding integral is estimated by
\begin{align}\label{Jm}
  \Big|\int_{0}^{\infty}\chi_{\infty}(|t|r)\,\chi_{\infty}(r)\,
  \widetilde{a}_{m}^{\pm}(r)\,e^{ir(t\pm s)}\,\df{r}\Big|\,
  \lesssim\,
  m^{\nu_{2}}\,\int_{\frac{1}{2}}^{\infty}
  r^{-2}\,\df{r}\,
  \lesssim\,
  m^{\nu_{2}}.
\end{align}

For the first term of the decomposition of $a_{m}^{\pm}(r)$, it is sufficient to obtain the estimate of the corresponding oscillatory integral $J_{0}(t,s)$ with amplitude $a_{0}^{\pm}(r)=(r^{2}+r_{0})^{-\theta/2}\mathbf{c}(\mp2r)^{-1}$, since the remaining part of $a_m^\pm(r)-\widetilde{a}_m^\pm = C_{n,m}a_0^\pm$. Note that $a_{0}^{\pm}(r)$ is a differential symbol of order $-1$ when $\re\theta=(n+1)/2$, we need more effort for this part.

Let us denote $\xi=t\pm s$ and divide our argument into two cases depending on whether $|\xi|$ is small or large. Assume at first that $|\xi|\ge1/2$. This case can be easily handled by using the integration by parts based on the fact that $i\xi e^{i\xi r}=\partial_{r}e^{i\xi r}$.  Formally, we have
\begin{align*}
  J_{0}(t,s)\,
  =\,\frac{i}{\xi}\,&\int_{0}^{\infty}e^{i\xi r}
  \,\frac{\partial}{\partial r}\,
  \Big(\chi_{\infty}(|t|r)\,\chi_{\infty}(r)\,
  a_{0}^{\pm}(r)\Big)\,\df{r}.
\end{align*}
According to \eqref{c estimate}, we have
\begin{align*}
  \Big|\frac{\partial}{\partial r}\,
  (a_{0}^{\pm}(r))\Big|\,
  \lesssim\,
  |\theta|\,r^{-\re\theta+\frac{n-1}{2}-1},
\end{align*}
for all $r\in\supp\chi_{\infty}$. Therefore, the contribution of the term where the derivative is applied to $a_{0}^{\pm}(r)$ is estimated by
\begin{align*}
  \frac{|\theta|}{|\xi|}\,
  \int_{\frac{1}{2}}^{\infty}r^{-2}\,\df{r}\,
  \lesssim\,|\theta|,
\end{align*}
provided that $|\xi|\ge1/2$ and $\theta\in\mathbb{C}$ with $\re{\theta}=\frac{n+1}{2}$.  If the derivative hits $\chi_{\infty}(r)$, then the corresponding integral is reduced to the interval $[1/2,1]$, which is clearly finite. In the remaining case where the derivative is applied to the cut-off function $\chi_{\infty}(|t|r)$, which only matters if $|t|<2$ according to \eqref{cut-off}, there is a constant $C>0$ such that the corresponding integral is estimated by 
\begin{align*}
  \frac{C|t|}{|\xi|}\,\int_{|t|^{-1}}^{2|t|^{-1}}r^{-1}\,\df{r}\,
  \le\,4(\ln2)C,
\end{align*}
provided that $|\xi|\ge1/2$ and $\theta\in\mathbb{C}$ with $\re{\theta}=\frac{n+1}{2}$. 

We conclude that in the case where $|\xi|\ge1/2$, we have the estimate
\begin{align}\label{large xi}
  |J_{0}(t,s)|\,
  \lesssim\,|\theta|
\end{align}
for $(t,s)\in\mathbb{R}^{*}\times\mathbb{R}_{+}^{*}$ satisfying $|t\pm s|\ge1/2$, and $\theta\in\mathbb{C}$ with $\re\theta=(n+1)/2$.

Now let us turn to the case where $|\xi|<1/2$. According to \eqref{c expansion} and the fact that 
\begin{align*}
  (r^{2}+r_{0})^{-\frac{\theta}{2}}\,
  =\,r^{-\theta}(1+\mathrm{O}(|\theta|\,r^{-2}))
  \qquad\forall\,r\in\supp\chi_{\infty},
\end{align*}
we know that, for $\theta\in\mathbb{C}$ with $\re\theta=(n+1)/2$, there exists a constant $C_{n}>0$ such that 
\begin{align*}
  a_{0}^{\pm}(r)\,
  =\,C_{n}\,r^{-1-i\im\theta}\,+\,R(r),
\end{align*}
with remainder $R(r)=\mathrm{O}(|\theta|\,r^{-2})$. Then we divide our integral into three parts $J_{0}=J_{0}^{-}+J_{0}^{+}+J_{0}^{R}$ with 
\begin{align*}
  J_{0}^{-}(t,s)\,+J_{0}^{+}(t,s)\,
  =\,C_{n}\,
  \Big(\int_{0}^{\frac{1}{|\xi|}}\,+\,\int_{\frac{1}{|\xi|}}^{\infty}\Big)
  \chi_{\infty}(|t|r)\,\chi_{\infty}(r)\,
  r^{-1-\im\theta}\,e^{i\xi r}\,\df{r},
\end{align*}
and 
\begin{align*}
  J_{0}^{R}(t,s)\,
  =\,
  \int_{0}^{\infty}\chi_{\infty}(|t|r)\,\chi_{\infty}(r)\,
  R(r)\,e^{i\xi r}\,\df{r}.
\end{align*}
Since $R(r)=\mathrm{O}(|\theta|\,r^{-2})$ for all $r\in\supp\chi_{\infty}$, it is easy to deduce that 
\begin{align}\label{JR}
  J_{0}^{R}(t,s)\,=\,\mathrm{O}(|\theta|).
\end{align}

Note that $1/|\xi|>2$ in the present case. To study $J_{0}^{-}$, we perform an integration by parts by using the formula $r^{-1-i\im\theta}=i(\im\theta)^{-1}\partial_{r}r^{-i\im\theta}$:
\begin{align*}
  J_{0}^{-}(t,s)\,
  &=\,C_{n}\,\int_{0}^{\frac{1}{|\xi|}}\,
  \chi_{\infty}(|t|r)\,\chi_{\infty}(r)\,
  r^{-1-i\im\theta}\,e^{i\xi r}\,\df{r}\\
  &=\,
  \frac{iC_{n}}{\im\theta}\,r^{-i\im\theta}\,
  \chi_{\infty}(|t|r)\,\chi_{\infty}(r)\,
  e^{i\xi r}\,\Big|_{r=0}^{r=|\xi|^{-1}}\\
  &\quad-\,
  \frac{iC_{n}}{\im\theta}\,
  \int_{0}^{\frac{1}{|\xi|}}\,r^{-i\im\theta}\,
  \frac{\partial}{\partial r}\,
  \Big(\chi_{\infty}(|t|r)\,\chi_{\infty}(r)\,
  e^{i\xi r}\Big)\,\df{r},
\end{align*}
where the first term is $\mathrm{O}(|\im\theta|^{-1})$. In the second integral, if the derivative is applied to the cut-off function $\chi_{\infty}(r)$, then the integral is reduced to the interval $[1/2,1]$, then bouned; while if the derivative hits $e^{i\xi r}$, then one can gain an extra $|\xi|$ to compensate for the $1/|\xi|$ from the upper bound of the integral. In the remaining case where the derivative is applied to the cut-off function $\chi_{\infty}(|t|r)$ when $|t|<2$, the corresponding term is estimated by 
  \begin{align*}
    C\,|t|\,\int_{|t|^{-1}}^{2|t|^{-1}}\df{r}\,=\,C,
  \end{align*}
for some $C>0$. We conclude that 
  \begin{align}\label{J0-}
    J_{0}^{-}(t,s)\,=\,\mathrm{O}(|\im\theta|^{-1}).
  \end{align}

  It remains for us to study $J_{0}^{+}$. Note that $\chi_{\infty}(r)=1$ in the present case where $r\ge1/|\xi|>2$.  By using the integration by parts based on the formula $e^{i\xi r}=(i\xi)^{-1}\partial_{r}e^{i\xi r}$ again, we obtain
  \begin{align*}
    J_{1}^{+}(t,s)\,
    &=\,C_{n}\,\int_{\frac{1}{|\xi|}}^{\infty}\,
    \chi_{\infty}(|t|r)\,r^{-1-i\im\theta}\,e^{i\xi r}\,\df{r}\\
    &=\,\frac{C_{n}}{i\xi}\,e^{i\xi r}\,
    r^{-1-i\im\theta}\,\chi_{\infty}(|t|r)\,
    \Big|_{r=|\xi|^{-1}}^{r=\infty}\\
    &\quad+\,\frac{iC_{n}}{\xi}\,
    \int_{\frac{1}{|\xi|}}^{\infty}\,e^{i\xi r}\,
    \frac{\partial}{\partial r}\,
    \Big(r^{-1-i\im\theta}\,
    \chi_{\infty}(|t|r)\Big),
  \end{align*}
  where the first term is $\mathrm{O}(1)$. The contribution of 
  \begin{align*}
    \frac{\partial}{\partial r}\,r^{-1-i\im\theta}\,
    =\,-(1+i\im\theta)\,r^{-2-i\im\theta}
  \end{align*}
  is also easy to handle. If the derivative is applied to the cut-off function $\chi_{\infty}(|t|r)$ when $|t|<2$, the support of $\partial_{r}\chi_{\infty}(|t|r)$ is reduced to the interval $[|t|^{-1},2|t|^{-1}]$, which means that the integral over $[|\xi|^{-1},\infty]$ vanishes unless $|t|/|\xi|<2$. In this case, the corresponding term is estimated by  
  \begin{align*}
    C\,\frac{|t|}{|\xi|}\,
    \int_{|t|^{-1}}^{2|t|^{-1}}\,r^{-1}\,\df{r}\,
    <\,2\ln2\,C,
  \end{align*}
  for some $C>0$. We conclude that 
  \begin{align}\label{J0+}
    J_{0}^{+}(t,s)\,=\,\mathrm{O}(1+|\im\theta|),
  \end{align}
  and deduce from \eqref{JR}, \eqref{J0-}, and \eqref{J0+} that 
  \begin{align}
    J_{0}(t,s)\,
    =\,\mathrm{O}\Big(|\im\theta|\,+\,\frac{1}{|\im\theta|}\Big),
  \end{align}
  for all $t\in\mathbb{R}^{*}$ and $s\in\mathbb{R}_{+}^{*}$ satisfying $|t\pm s|<1/2$, and $\theta\in\mathbb{C}$ with $\re\theta=(n+1)/2$. Together with the estimate \eqref{large xi} for the case $|t\pm s|\ge1/2$, and the estimate \eqref{Jm} for $m\ge1$, we finally obtain 
  \begin{align*}
    |J_{m}(t,s)|\,
    \lesssim\,(1+m^{\nu_{1}})\Big(|\im\theta|\,+\,\frac{1}{|\im\theta|}\Big)\,
    +\,m^{\nu_{2}}
  \end{align*}
  for all $m\in\mathbb{N}$. The lemma is therefore proved by taking $\nu=\max\lbrace\nu_{1},\nu_{2}\rbrace$.
\end{proof}

\begin{proposition}\label{estimate Iinf+odd}
  Let $n$ be odd and $\theta\in\mathbb{C}$ with $\re{\theta}=\frac{n+1}{2}$. Then, for all $x\in G$ and $t\in\mathbb{R}^{*}$ such that $|x|/|t|>1/2$, there exists $C_{n,\theta}>0$ such that 
  \begin{align*}
    |I_{+}^{\infty}(t,x)|\,
    \le\,C_{n,\theta}\,e^{-\frac{n-1}{2}|x|}.
  \end{align*}
\end{proposition}

\begin{proof}
We will sketch the proof to avoid repeating the same calculations done in even dimensions. In the odd dimension, the expressions for the scalar components $\varphi_{\pm}^{+}(r,\cdot)$ and $\varphi_{\pm}^{-}(r,\cdot)$ of the spherical function $\Phi_{\pm}(r,\cdot)$ are more complicated. They are given in \cite[Theorem 5.4]{CP01}:
\begin{align*}
  \varphi_{\pm}^{+}(r,s)\,
  &=\,\Big(\cosh\frac{s}{2}\Big)\,
  \phi_{2r}^{(\frac{n}{2}-1,\frac{n}{2})}\Big(\frac{s}{2}\Big)\,
  \mp\,\frac{2ri}{n}\,\Big(\sinh\frac{s}{2}\Big)\,
  \phi_{2r}^{(\frac{n}{2},\frac{n}{2}-1)}\Big(\frac{s}{2}\Big),\\
  \varphi_{\pm}^{-}(r,s)\,
  &=\,\Big(\cosh\frac{s}{2}\Big)\,
  \phi_{2r}^{(\frac{n}{2}-1,\frac{n}{2})}\Big(\frac{s}{2}\Big)\,
  \pm\,\frac{2ri}{n}\,\Big(\sinh\frac{s}{2}\Big)\,
  \phi_{2r}^{(\frac{n}{2},\frac{n}{2}-1)}\Big(\frac{s}{2}\Big).
\end{align*}
The first part of $\varphi_{\pm}^{+}$ and $\varphi_{\pm}^{-}$, which is identical to $\varphi^{\pm}$ defined in \eqref{scalar even}, can be treated in the same way as in the case of the even dimension. Let us denote by
\begin{align*}
  \widetilde{\varphi}(r,s)\,
  =\,\frac{2ri}{n}\,\Big(\sinh\frac{s}{2}\Big)\,
  \phi_{2r}^{(\frac{n}{2},\frac{n}{2}-1)}\Big(\frac{s}{2}\Big)
\end{align*}
the second part of $\varphi_{\pm}^{+}$ and $\varphi_{\pm}^{-}$. According to the proof of \cite[Theorem 5.4]{CP01}, we know that
\begin{align*}
  ir\,\widetilde{\varphi}(r,s)\,
  =\,\sqrt{-z}\,
  \Big(\frac{n}{2}-(1-z)\frac{\partial}{\partial z}\Big)\,
  {}_{2}F_{1}\Big(\frac{n}{2}+ir,\frac{n}{2}-ir;\frac{n}{2};z\Big)
\end{align*}
with $z=-(\sinh(s/2))^{2}$. Note that, for all $s>0$,
\begin{align*}
  (1-z)\frac{\partial}{\partial z}\,
  =\,\Big(-\coth\frac{s}{2}\Big)\,\frac{\partial}{\partial s}.
\end{align*}
Together with \eqref{scalar even}, we get 
\begin{align*}
  \widetilde{\varphi}(r,s)\,
  &=\,\frac{1}{ir}\,
  \Big(\sinh\frac{s}{2}\Big)\,
  \Big(\frac{n}{2}\,+\,
  \Big(\coth\frac{s}{2}\Big)\,\frac{\partial}{\partial s}\Big)\,
  \phi_{2r}^{(\frac{n}{2}-1,\frac{n}{2})}\Big(\frac{s}{2}\Big)\\
  &=\,\frac{n}{2ir}\,\Big(\sinh\frac{s}{2}\Big)\,
  \phi_{2r}^{(\frac{n}{2}-1,\frac{n}{2})}\Big(\frac{s}{2}\Big)\,
  +\,\frac{1}{ir}\,\Big(\cosh\frac{s}{2}\Big)\,
  \frac{\partial}{\partial s}\,
  \phi_{2r}^{(\frac{n}{2}-1,\frac{n}{2})}\Big(\frac{s}{2}\Big),
\end{align*}
where the first term can be estimated again using \eqref{HC expansion1} and \eqref{HC expansion2}. Note that the hyperbolic sine function also satisfies $\sinh(s/2)\lesssim\,e^{s/2}$ for $s>0$. Then we return to the analysis in even dimension. In fact, the analysis is even simpler here because of the factor $r^{-1}$, which provides more decay for all $r\in\supp\chi_{\infty}$. 

Regarding the second term of $\widetilde{\varphi}(r,s)$, we note that the derivative is actually applied to the functions $\Phi_{r}(s)$ and $\Phi_{-r}(s)$:
\begin{align*}
  \frac{\partial}{\partial s}\,
  \phi_{2r}^{(\frac{n}{2}-1,\frac{n}{2})}\Big(\frac{s}{2}\Big)\,
  =\,
  \mathbf{c}(2r)\,\frac{\partial}{\partial s}\,
  \Phi_{2r}\Big(\frac{s}{2}\Big)\,
  +\,
  \mathbf{c}(-2r)\,\frac{\partial}{\partial s}\,
  \Phi_{2r}\Big(\frac{s}{2}\Big),
\end{align*}
according to the expansion \eqref{HC expansion1}. We know from \eqref{HC expansion2} that the derivative gives either $-m/2$ in the series or the term $ir-n/2$. On the one hand, any additional polynomial factor in the series defined in \eqref{HC expansion2} has no essential contribution. On the other hand, the contribution of $ir-n/2$ is compensated by the factor $(ir)^{-1}$. The proof then follows with the same calculations in even dimension.
\end{proof}

Let us divide the kernel $K_{t}^{\theta}=K_{t,0}^{\theta}+K_{t,\infty}^{\theta}$ into two parts with
\begin{align*}
  K_{t,\,0}^{\theta}\,=\,\sum_{\sigma\in\widehat{M}(\tau)}I_{0}(t,x)
  \qquad\textnormal{and}\qquad
  K_{t,\,\infty}^{\theta}\,=\,\sum_{\sigma\in\widehat{M}(\tau)}I_{\infty}(t,x),
\end{align*}
where the integrals $I_{0}(t,x)$ and $I_{\infty}(t,x)$ also depend on the parameter $\theta$ and the representation $\sigma$. Note that the estimates of $I_{\infty}(t,x)$ obtained for $\re\theta=(n+1)/2$ also hold in the simpler case where $\re\theta>(n+1)/2$. According to the Propositions \ref{estimate I0}, \ref{estimate Iinf-}, \ref{estimate Iinf+-}, \ref{estimate Iinf+even}, and \ref{estimate Iinf+odd} proved above, we summarize the following corollary.

\begin{corollary}\label{kernel corollary}
  Let $\tau$ be either the spin representation $\tau_{n}$ of $K$ on $V_{\tau_{n}}$ in odd dimension, or the half-spin representations $\tau_{n}^{\pm}$ of $K$ on $V_{\tau_{n}^{\pm}}$ in even dimension. Suppose that $\theta\in\mathbb{C}$. Then there exists $C_{n,\theta}>0$ such that $C_{n,\theta}=\mathrm{O}(|\im\theta|+|\im\theta|^{-1})$ and the following pointwise estimate hold for all $x\in G$ and $t\in\mathbb{R}^{*}$:
  \begin{align*}
    |K_{t,\,0}^{\theta}(x)|\,
    \le\,
    C_{n,\theta}\,\varphi_{0}(x)\,
    \begin{cases}
      1
      &\qquad\textnormal{if}\quad |t|<1,\\
      |t|^{-1}\,(1+|x|)
      &\qquad\textnormal{if}\quad |t|\ge1.
    \end{cases}
  \end{align*}
  Assume in addition that $\theta\in\mathbb{C}$ with $\re\theta\ge(n+1)/2$. Then
  \begin{align*}
    |K_{t,\,\infty}^{\theta}(x)|\,
    \le\,
    C_{n,\theta}\,e^{-\frac{n-1}{2}|x|}\,
    \begin{cases}
      |t|^{-\frac{n-1}{2}}
      &\qquad\textnormal{if}\quad |t|<1,\\
      |t|^{-N}\,(1+|x|)^{N}
      &\qquad\textnormal{if}\quad |t|\ge1,
    \end{cases}
  \end{align*}
  for any $N\ge0$ and for all $x\in G$ and $t\in\mathbb{R}^{*}$.
\end{corollary}

\section{Dispersive property and Strichartz inequality}\label{sec: strichartz}
In this section, we analyze the dispersive properties of the Dirac propagator using the kernel estimates derived previously and establish the Strichartz inequality  for \eqref{Dirac0} through the $TT^{*}$ argument.

The approach involves interpolating the analytic family of smoothed Dirac propagators $\mathbf{D}_{t}^{\theta}=(D^{2}-\mathcal{R}/4)^{-\theta/2}e^{it|D|}$ within the vertical strip $0\le\re\theta\le(n+1)/2$. However, as noted in Corollary \ref{kernel corollary}, the constant $C_{n,\theta}$ in the pointwise kernel estimates has the singularity $|\im\theta|^{-1}$. To deal with this, the standard way is to introduce an additional constant factor
\begin{align*}
  \widetilde{C}_{n,\theta}\,
  :=\,\frac{e^{\theta^{2}}}{\Gamma(\frac{n+1}{2}-\theta)},
\end{align*}
where the Gamma function helps cancel such a singularity  at the boundary point  $\re\theta=(n+1)/2$, while the exponential function ensures boundedness at infinity within the vertical strip. More precisely, we have 
\begin{align*}
  |\widetilde{C}_{n,\theta}|\,
  \lesssim\,
  \Big|\theta-\frac{n+1}{2}\Big|\,|\theta|\,
  e^{\pi|\im\theta|-|\im\theta|^{2}},
\end{align*}
see, for instance, \cite[(3.2)]{AZ24}. Therefore, we will give the dispersive estimates for the analytic family of operators $\lbrace\widetilde{\mathbf{D}}_{t}^{\theta}\rbrace_{\theta\in\mathbb{C}}$ where
\begin{align*}
  \widetilde{\mathbf{D}}_{t}^{\theta}\,
  \,=\,
  \widetilde{C}_{n,\theta}\,\mathbf{D}_{t}^{\theta}\,
  =\,
  \frac{e^{\theta^{2}}}{\Gamma(\frac{n+1}{2}-\sigma)}\,
  \left(D^{2}-\frac{\mathcal{R}}{4}\right)^{-\frac{\theta}{2}}\,e^{it|D|},
\end{align*}
in the vertical strip $0\le\re\theta\le(n+1)/2$. Recall that $\mathcal{R}=-n(n-1)$ denotes the normalized scalar curvature of $\mathbb{H}^{n}$.

Another key ingredient is an available norm estimate for the convolution product over spinor bundles. Let $F\in\mathcal{S}(G,\tau)$ and $\kappa\in\mathcal{S}(G,\tau,\tau)$, the convolution product $F*\kappa\in\mathcal{S}(G,\tau)$ on spinor bundles is defined by
\begin{align*}
  (F*\kappa)(x)\,
  =\,
  \int_{G}\kappa(y^{-1}x)\,F(y)\,\df{y}.
\end{align*}
The following proposition gives the $L^{q'}$-$L^{q}$ norm estimate for the convolution product, and can be seen as the vector-valued Kunze-Stein phenomenon. Its proof is based on the last theorem of \cite{Ste70}.

\begin{proposition}\label{KS}
  Let $\mathcal{A}_{\infty}(G,\tau,\tau)$ be the space of $\tau$-radial functions in $L^{\infty}(G,\tau,\tau)$ equipped with the $L^{\infty}$ norm, and $\mathcal{A}_{q}(G,\tau,\tau)$ be the space of $\tau$-radial functions in $\varphi_{0}$-weighted $L^{q/2}(G,\tau,\tau)$ equipped with the norm 
  \begin{align*}
    \|\cdot\|_{\mathcal{A}_{q}}\,
    =\,
    \Big(
      \int_{G}|\cdot|^{\frac{q}{2}}\,\varphi_{0}(x)\,\df{x}
    \Big)^{\frac{2}{q}}
    \qquad\forall\,2\le q<\infty.
  \end{align*}
  Then the inclusion
  \begin{align}\label{KS inclusion}
    L^{q'}(G,\tau)\,*\,\mathcal{A}_{q}(G,\tau,\tau)\,\subset\,L^{q}(G,\tau)
  \end{align}
  holds for all $2\le q\le\infty$.
\end{proposition}

\begin{proof}
  The inclusion \eqref{KS inclusion} simply holds for $q=\infty$, since
  \begin{align*}
    |(F*\kappa)(x)|\,
    \le\,
    \int_{G}|\kappa(y^{-1}x)|\,|F(y)|\,\df{y}\,
    \le\,
    \|\kappa\|_{\mathcal{A}_{\infty}}\,
    \|F\|_{L^{1}(G,\tau)},
  \end{align*}
  for all $F\in L^{1}(G,\tau)$. In the case where $q=2$, we consider the operator $T_{\kappa}:\mathcal{S}(G,\tau)\rightarrow\mathcal{S}(G,\tau)$ defined by $T_{\kappa}F=F*\kappa$, where $\kappa$ is a $\tau$-radial function on $G$. Recall that the scalar ground spherical function $\varphi_{0}$ defined in \eqref{phi0} is positive and bi-$K$-invariant. We write
  \begin{align*}
    T_{\kappa}(\varphi_{0}F)(x)\,
    &=\,\int_{G}\kappa(y^{-1}x)\,(\varphi_{0}F)(y)\,\df{y}\\
    &=\,\int_{G}\kappa(zk^{-1})\,(\varphi_{0}F)(xkz^{-1})\,\df{z}
  \end{align*}
  by substituting $y=xkz^{-1}$ with $k\in K$. Then 
  \begin{align}\label{Tkappa}
    |T_{\kappa}(\varphi_{0}F)(x)|\,
    &\le\,\int_{G}
    |\tau(k^{-1})^{-1}\kappa(z)|\,
    |(\varphi_{0}F)(xkz^{-1})|\,\df{z}\notag\\
    &\lesssim\,
    \Big(
      \int_{G}
      |\kappa(z)|\,\varphi_{0}(xkz^{-1})\,\df{z}
    \Big)\,
    \|F\|_{L^{\infty}(G,\tau)},
  \end{align}
  since $\kappa$ is $\tau$-radial, where $\tau$ is a bounded representation of $K$. Recall that for any $r>0$ the spherical function $\varphi_{r}$ satisfies
  \begin{align*}
    \int_{K}\varphi_{r}(xkz^{-1})\,\df{k}\,
    =\,
    \varphi_{r}(x)\,\varphi_{-r}(z)
    \qquad\forall\,x,z\in G.
  \end{align*}
  We deduce by integrating $k$ over $K$ on both sides of \eqref{Tkappa} that 
  \begin{align*}
    |T_{\kappa}(\varphi_{0}F)(x)|\,
    \lesssim\,
    \varphi_{0}(x)\,
    \Big(
      \int_{G}
      |\kappa(z)|\,\varphi_{0}(z)\,\df{z}
    \Big)\,
    \|F\|_{L^{\infty}(G,\tau)}.
  \end{align*}
  In other words, we obtain
  \begin{align*}
    \|\varphi_{0}^{-1}T_{\kappa}F\|_{L^{\infty}(G,\tau)}\,
    \lesssim\,
    \|\kappa\|_{\mathcal{A}_{2}(G,\tau,\tau)}\,
    \|\varphi_{0}^{-1}F\|_{L^{\infty}(G,\tau)},
  \end{align*}
  and 
  \begin{align*}
    \|\varphi_{0}T_{\kappa}F\|_{L^{1}(G,\tau)}\,
    \lesssim\,
    \|\kappa\|_{\mathcal{A}_{2}(G,\tau,\tau)}\,
    \|\varphi_{0}F\|_{L^{1}(G,\tau)},
  \end{align*}
  by duality. By interpolating between the last two inequalities, we deduce that the inclusion \eqref{KS inclusion} holds for $q=2$. The proposition follows from the interpolation between
  \begin{align*}
    L^{2}(G,\tau)\,*\,\mathcal{A}_{2}(G,\tau,\tau)\,\subset\,L^{2}(G,\tau),
  \end{align*}
  and 
  \begin{align*}
    L^{1}(G,\tau)\,*\,\mathcal{A}_{\infty}(G,\tau,\tau)\,\subset\,L^{\infty}(G,\tau),
  \end{align*}
  by noting that the interpolation space between $\mathcal{A}_{2}$ and $\mathcal{A}_{\infty}$ is $\varphi_{0}$-weighted $L^{q/2}(G,\tau)$ space, i.e., $\mathcal{A}_{q}(G,\tau,\tau)$.
\end{proof}

Once again, let $\tau$ be either the spin representation $\tau_{n}$ of $K$ on $V_{\tau_{n}}$ in odd dimension, or the half-spin representations $\tau_{n}^{\pm}$ of $K$ on $V_{\tau_{n}^{\pm}}$ in even dimension. Now we prove the dispersive estimate for the operator $\widetilde{\mathbf{D}}_{t}^{\theta}$. Theorem \ref{main thm dispersive} follows from the following result.

\begin{theorem}\label{thm: dispersive}
  Let $2<q<\infty$ and $\theta\ge(n+1)(1/2-1/q)$. Then the operator $\widetilde{\mathbf{D}}_{t}^{\theta}$ is bounded from $L^{q'}(G,\tau)$ to $L^{q}(G,\tau)$:
  \begin{align*}
    \|\widetilde{\mathbf{D}}_{t}^{\theta}
    \|_{L^{q'}(G,\tau)\rightarrow L^{q}(G,\tau)}\,
    \lesssim\,
    \begin{cases}
      |t|^{-(n-1)(\frac{1}{2}-\frac{1}{q})}
      &\qquad\textnormal{if}\quad |t|<1,\\
      |t|^{-1}
      &\qquad\textnormal{if}\quad |t|\ge1.
    \end{cases}
  \end{align*}
\end{theorem}

\begin{proof}
  For simplicity, we restrict to the critical case where $\theta=(n+1)(1/2-1/q)$. We divide the proof into two parts according to the decomposition $\widetilde{\mathbf{D}}_{t}^{\theta}=\widetilde{\mathbf{D}}_{t,0}^{\theta}+\widetilde{\mathbf{D}}_{t,\infty}^{\theta}$, where $\widetilde{\mathbf{D}}_{t,0}^{\theta}$ and $\widetilde{\mathbf{D}}_{t,\infty}^{\theta}$ are operators associated with the convolution kernels $K_{t,\,0}^{\theta}$ and $K_{t,\,\infty}^{\theta}$ in Corollary \ref{kernel corollary}. For $\widetilde{\mathbf{D}}_{t,0}^{\theta}$, we derive the results by direct computations based on the Corollary \ref{kernel corollary}, while for $\widetilde{\mathbf{D}}_{t,\infty}^{\theta}$, we use an analytic interpolation argument in the vertical strip $0\le\re\theta\le(n+1)/2$.

  Since $K_{t,\,0}^{\theta}$ is a $\tau$-radial function, using Proposition \ref{KS} and Corollary \ref{kernel corollary} we get
  \begin{align*}
    \|\widetilde{\mathbf{D}}_{t,0}^{\theta}
    \|_{L^{q'}(G,\tau)\rightarrow L^{q}(G,\tau)}\,
    &\lesssim\,
    \Big(
      \int_{G}|\,K_{t,\,0}^{\theta}(x)|^{\frac{q}{2}}\,\varphi_{0}(x)\,\df{x}
    \Big)^{\frac{2}{q}}\\
    &\lesssim\,
    \min\lbrace{1,|t|^{-1}}\rbrace\,
    \Big(
      \int_{G}\,(1+|x|)^{\frac{q}{2}}
      \varphi_{0}(x)^{\frac{q}{2}+1}\,\df{x}
    \Big)^{\frac{2}{q}},
  \end{align*}
  for all $\theta\in\mathbb{C}$. According to the Cartan decomposition \eqref{CartanHaar} and the estimate \eqref{phi0} of the ground spherical function, we obtain 
  \begin{align*}
    \|\widetilde{\mathbf{D}}_{t,0}^{\theta}
    \|_{L^{q'}(G,\tau)\rightarrow L^{q}(G,\tau)}\,
    \lesssim\,
    \min\lbrace{1,|t|^{-1}}\rbrace
    \Big(
      \int_{0}^{\infty}(1+r)^{q+1}\,
      e^{-\frac{n-1}{2}(\frac{q}{2}-1)r}\,\df{r}
    \Big)^{\frac{2}{q}},
  \end{align*}
  which is finite provided that $q>2$. 

  For the operator $\widetilde{\mathbf{D}}_{t,\infty}^{\theta}$ we know from the spectral theorem that, if $\re\theta=0$,
  \begin{align*}
    \|\widetilde{\mathbf{D}}_{t,\infty}^{\theta}
    \|_{L^{2}(G,\tau)\rightarrow L^{2}(G,\tau)}
    \lesssim\,1.
  \end{align*}
  On the other hand, we deduce from Corollary \ref{kernel corollary} that, if $\re\theta=\frac{n+1}{2}$,
  \begin{align*}
    \|\widetilde{\mathbf{D}}_{t,\infty}^{\theta}
    \|_{L^{1}(G,\tau)\rightarrow L^{\infty}(G,\tau)}
    &\lesssim\,
    \|K_{t,\,\infty}^{\theta}
    \|_{L^{\infty}(G,\tau)}\\
    &\lesssim\,
    \begin{cases}
      |t|^{-\frac{n-1}{2}}
      &\qquad\textnormal{if}\quad |t|<1,\\
      |t|^{-N}
      &\qquad\textnormal{if}\quad |t|\ge1,
    \end{cases}
  \end{align*}
  for any $N\ge0$. By applying the interpolation theorem to the analytic family of operators $\widetilde{\mathbf{D}}_{t}^{\theta}$ in the vertical strip $0\le\re\theta\le(n+1)/2$, and by combining the estimate concerning $K_{t,\,0}^{\theta}$, we prove the dispersive estimate
  \begin{align*}
    \Big\|\widetilde{\mathbf{D}}_{t}^{\frac{n+1}{2}(1-\frac{2}{q})}
    \Big\|_{L^{q'}(G,\tau)\rightarrow L^{q}(G,\tau)}\,
    \lesssim\,
    \begin{cases}
      1\,+\,|t|^{-\frac{n-1}{2}(1-\frac{2}{q})}
      &\qquad\textnormal{if}\quad |t|<1,\\
      |t|^{-1}\,+\,|t|^{-N(1-\frac{2}{q})}
      &\qquad\textnormal{if}\quad |t|\ge1,
    \end{cases}
  \end{align*}
  for any $N\ge0$. The theorem is therefore proved.
\end{proof}

Let $s\in\mathbb{R}$ and $1<q<+\infty$. We define the Sobolev space $H^{s,q}(G,\tau)$ as the image of $L^{q}(G,\tau)$ under the operator $(D^{2}-\mathcal{R}/4)^{-s/2}$, equipped with the norm
\begin{align*}
  \|F\|_{H^{s,q}(G,\tau)}\,
  :=\,
  \left\|\left(D^{2}-\frac{\mathcal{R}}{4}
  \right)^{s/2}F\right\|_{L^{q}(G,\tau)}.
\end{align*}
Recall that $\mathcal{R}=-n(n-1)$ denotes the normalized scalar curvature of $\mathbb{H}^{n}$. We will use the following Sobolev embedding result.

\begin{proposition}\label{prop: sobolev}
  Let $1<q_1 \le q_2<+\infty$ and $s_1 \ge s_2\ge0$ such that $s_1-s_2=n(1/q_1 - 1/q_2)$. Then
  \begin{align*}
    H^{s_{1},q_{1}}(G,\tau)\,
    \subset\,H^{s_{2},q_{2}}(G,\tau).
  \end{align*}
\end{proposition}

\begin{proof}
  According to the Lichnerowicz formula, we have
  \begin{align*}
    D^{2}\,=\,\nabla^{*}\nabla + \frac{\mathcal{R}}{4},
  \end{align*}
  where $\nabla^{*}\nabla$ represents the Bochner Laplacian acting on spinors. It follows from \cite[Section 6]{Str83} that $\|(\nabla^{*}\nabla)^{1/2}F\|_{L^{q}}$ and $\|\nabla F\|_{L^{q}}$ are equivalent for all $1<q<\infty$ on non-comapct symmetric spaces of rank $1$, which include the real hyperbolic space $\mathbb{H}^{n}$. We then have the equivalence between the Sobolev norms $\|F\|_{H^{s,q}(G,\tau)}$, defined by the Dirac operator, and 
  \begin{align*}
    \|F\|_{W^{s,q}(G,\tau)}\,
    =\,
    \left( \sum_{j=0}^s 
    \int_G |\nabla^j F(g)|^q \, 
    \df{g} \right)^{1/q},
  \end{align*}
  defined by the covariant derivative, for $s=1$. The equivalence for $s>1$ is more complicated due to the non-commutativity of the operators $\nabla^{*}\nabla$ and $\nabla$. For the scalar case where $\nabla^{*}\nabla$ is the Laplace-Beltrami operator, we refer to \cite[Section 7.4.5]{Tri92}, and for the case of vector bundles, to \cite{Yos92}, since $\mathbb{H}^n$ has a bounded geometry.
  
 According to \cite[Lemma 3.1 and Theorem 3.14]{Heb96}, the Sobolev embedding
  \begin{align*}
    W^{s_{1},q_{1}}(G,\tau)\,
    \subset\,
    W^{s_{2},q_{2}}(G,\tau)
  \end{align*}
  holds for all $1/q_2=1/q_1 - (s_1-s_2)/n$ with $s_1$, $s_2$ nonnegative integers. The proposition then follows by interpolation.
\end{proof}

Let us now turn to the Strichartz inequality. The solution to the linear Dirac equation
\begin{align}\label{DiracLinear}
  \begin{cases}
    (i\partial_{t}\,+\,D)\,U(t,x)\,=\,L(t,x),\\
    U(0,x)\,=\,U_{0}(x),
  \end{cases}
\end{align}
is given by 
\begin{align*}
  U(t,x)\,
  =\,
  e^{itD}U_{0}(x)\,-\,
  i\int_{0}^{t} e^{i(t-s)D} L(s,x)\,\df{s}.
\end{align*}
This solution is estimated, in a time-space mixed norm, by the initial data $U_{0}$ and $L$. Such an estimate is known as the Strichartz inequality and is a key tool in the study of non-linear problems.

\begin{theorem}\label{thm: strichartz}
  Let $(p_1,q_1)$ and $(p_2,q_2)$ be two admissible pairs, i.e., the points  $(1/p_{1},1/q_{1})$ and $(1/p_{2},1/q_{2})$ belong to the admissible set
  \begin{align}\label{def: admissible}
    \left\lbrace
    \left(\frac{1}{p},\frac{1}{q}\right)\in
    \left(0,\frac{1}{2}\right)\times
    \left(0,\frac{1}{2}\right)\,
    \Big|\,
    \frac{1}{p}\ge \frac{n-1}{2}
    \left(\frac{1}{2}-\frac{1}{q}\right)
    \right\rbrace
    \bigcup
    \left\lbrace
    \left(0,\frac{1}{2}\right)
    \right\rbrace.
  \end{align}
  Let 
  \begin{align*}
    \theta_{1}\,\ge\,\frac{n+1}{2}\left(\frac{1}{2}-\frac{1}{q_1}\right)
    \quad\textnormal{and}\quad
    \theta_{2}\,\ge\,\frac{n+1}{2}\left(\frac{1}{2}-\frac{1}{q_2}\right).
  \end{align*}
  Then, the solution of the linear Dirac equation \eqref{DiracLinear} satisfies the inequality
  \begin{align}
    \|U\|_{L^{p_{1}}(\mathbb{R};\,H^{-\theta_{1},q_{1}}(G,\tau))}\,
    \lesssim\,
    \|U_0\|_{L^{2}(G,\tau)}\,+\,
    \|L\|_{L^{p_{2}'}(\mathbb{R};\,H^{\theta_{2},q_{2}'}(G,\tau))}.
  \end{align}
\end{theorem}

\begin{proof}
  Once again, we denote $r_{0}=-\mathcal{R}/4$  for simplicity. Let 
  \begin{align*}
    TF(t,x)\,=\,
    (D^2+r_{0})^{-\frac{\theta}{2}} e^{\pm it|D|} F(x),
  \end{align*}
  and denote its adjoint operator by $T^*$, which is given by 
  \begin{align*}
    T^* F(x)\,=\,
    \int_{\mathbb{R}}\,
    (D^2+r_{0})^{-\frac{\theta}{2}} 
    e^{\mp is|D|} F(s,x)\,\df{s}.
  \end{align*}

  We want to prove that the operator $TT^{*}$ defined by
  \begin{align*}
    TT^{*}F(t,x)\,
    :=\,
    \int_{\mathbb{R}}\,
    (D^2+r_{0})^{-\theta} 
    e^{i(t\pm s)|D|} F(s,x)\,\df{s},
  \end{align*}
  is bounded from $L^{p'}(\mathbb{R};\,L^{q'}(G,\tau))$ to $L^{p}(\mathbb{R};\,L^{q}(G,\tau))$, for some pairs $(p,q)$ admissible. This is clearly true for $(p,q)=(+\infty,2)$ by the $L^{2}$ conservation.

  Together with the Minkowski integral inequality, we know from Theorem \ref{thm: dispersive} that $\|TT^* F\|_{L^{p}(\mathbb{R};\,L^{q}(G,\tau))}$ is bounded above by the sum of 
  \begin{align}\label{inproof: strichartz small}
    \left\| \int_{|t\pm s|<1} 
    |t\pm s|^{-(n-1)(\frac{1}{2}-\frac{1}{q})}\,
    \|F(s,\cdot)\|_{L^{q'}(G,\tau)}\,\df{s}
    \right\|_{L^{p}(\mathbb{R})}
  \end{align}
  and 
  \begin{align}\label{inproof: strichartz large}
    \left\| \int_{|t\pm s|\ge1} 
    |t\pm s|^{-1}\,
    \|F(s,\cdot)\|_{L^{q'}(G,\tau)}\,\df{s}
    \right\|_{L^{p}(\mathbb{R})},
  \end{align}
  up to a constant. On the one hand, the convolution kernel $|x|^{-1}\mathbf{1}_{\lbrace |x|\ge1 \rbrace}$ is in $L^{p/2}(\mathbb{R})$ for all $p>2$, then it defines a bounded operator from $L^{p'}(\mathbb{R})$ to $L^{p}(\mathbb{R})$ provided that $p>2$. On the other hand, the convolution kernel $|x|^{-(n-1)(\frac{1}{2}-\frac{1}{q})}\mathbf{1}_{\lbrace |x|<1 \rbrace}$ defines a bounded operator from $L^{p'}(\mathbb{R})$ to $L^{p}(\mathbb{R})$ provided that
  \begin{align*}
    0\,\le\,\frac{1}{p'}-\frac{1}{p}\,\le\,
    1-(n-1)\left(\frac{1}{2}-\frac{1}{q}\right).
  \end{align*}

  It thus follows that the operator $TT^{*}$ is bounded from $L^{p'}(\mathbb{R};\,L^{q'}(G,\tau))$ to $L^{p}(\mathbb{R};\,L^{q}(G,\tau))$ for all $(p,q)$ admissible. This, in turn, implies that the operator $T$ is bounded from $L^2(G,\tau)$ to $L^{p}(\mathbb{R};\,L^{q}(G,\tau))$, and the operator $T^*$ is bounded from $L^{p'}(\mathbb{R};\,L^{q'}(G,\tau))$ to $L^2(G,\tau)$.

  We deduce that, for $(p_1,q_1)$ admissible and $\theta_1\ge(n/2+1/2)(1/2-1/q_1)$,
 \begin{align*}
\|U\|_{L^{p_{1}}(\mathbb{R};\,H^{-\theta_{1},q_{1}}(G,\tau))}\,
    &=\,
    \|(D^2+r_0)^{-\frac{\theta_{1}}{2}}e^{itD}U_{0}
    \|_{L^{p_{1}}(\mathbb{R};\,L^{q_{1}}(G,\tau))}\\
    &\leq \,
    \|(D^2+r_0)^{-\frac{\theta_{1}}{2}}e^{it|D|}U_{0}^+
    \|_{L^{p_{1}}(\mathbb{R};\,L^{q_{1}}(G,\tau))}\\
    &\quad+  \|(D^2+r_0)^{-\frac{\theta_{1}}{2}}e^{-it|D|}U_{0}^-
    \|_{L^{p_{1}}(\mathbb{R};\,L^{q_{1}}(G,\tau))}\\
    &\lesssim \|U_{0}^+\|_{L^{2}(G,\tau)}+ \|U_{0}^-\|_{L^{2}(G,\tau)}\lesssim  \|U_{0}\|_{L^{2}(G,\tau)},
  \end{align*}
where $U_0^\pm =P^\pm U_0$ and $P^\pm$, defined by \eqref{def: proj operator}, are bounded operators on $L^2(G,\tau)$. This confirms the Strichartz inequality for the homogeneous part.

The inhomogeneous part involves studying the truncated $TT^*$ operator by combining the argument above with the Christ-Kiselev argument outlined in \cite{CK01}. This parallels the Euclidean case and is therefore omitted. In fact, the reasoning here is simpler, as the endpoint $p=2$, which presents technical difficulties, is already excluded.
\end{proof}

Using a standard argument based on the Sobolev embedding, the admissible set \eqref{def: admissible} can be extended to the entire square (see Figure \ref{fig: admissible}), where the Strichartz estimate remains valid. Theorem \ref{main thm strichartz} follows from the following corollary.

\begin{corollary}\label{cor: strichartz}
  Let $(p_1,q_1)$ and $(p_2,q_2)$ be two points in the admissible square 
  \begin{align}\label{def: admissible square}
    \left[0,\frac{1}{2}\right)\,
    \times\,
    \left(0,\frac{1}{2}\right)\,
    \bigcup\,
    \left\lbrace
    \left(0,\frac{1}{2}\right)
    \right\rbrace.
  \end{align}
  Denote by 
  \begin{align*}
    \theta_i(p_i, q_i)\,
    &=\,
    \frac{n+1}{2} \left( \frac{1}{2} - \frac{1}{q_{i}} \right)\,
    +\,
    \max \left\{ 0, \frac{n-1}{2} 
    \left( \frac{1}{2} - \frac{1}{q_i} \right) 
      - \frac{1}{p_i} \right\} \\
    &= 
  \begin{cases} 
    \frac{n+1}{2} 
    \left( \frac{1}{2} - \frac{1}{q_i} \right)
    & \text{if } \frac{1}{p_i} \geq 
    \frac{n-1}{2} \left( \frac{1}{2} - \frac{1}{q_i} \right), \\ 
    n \left( \frac{1}{2} - \frac{1}{q_i} \right) - \frac{1}{p_i}
    & \text{if } \frac{1}{p_i} \leq 
    \frac{n-1}{2} \left( \frac{1}{2} - \frac{1}{q_i} \right), 
  \end{cases}
  \end{align*}
  for $i=1,2$, and let $\theta_i\ge\theta_i(p_i, q_i)$. Then, the Strichartz inequality still holds, i.e., the solution of the linear Dirac equation \eqref{DiracLinear} satisfies the inequality
  \begin{align*}
    \|U\|_{L^{p_{1}}(\mathbb{R};\,H^{-\theta_{1},q_{1}}(G,\tau))}\,
    \lesssim\,
    \|U_0\|_{L^{2}(G,\tau)}\,+\,
    \|L\|_{L^{p_{2}'}(\mathbb{R};\,H^{\theta_{2},q_{2}'}(G,\tau))}.
  \end{align*}
\end{corollary}

\begin{proof}
This corollary is derived by combining the Strichartz inequality in Theorem \ref{thm: strichartz} with the Sobolev embedding in Proposition \ref{prop: sobolev}, and its proof can be found in many references. For the sake of completeness, we provide a brief sketch of the proof below.

For simplicity, we focus on the case where $\theta_i=\theta_i(p_i,q_i)$ for $i=1,2$. We choose two parameters $Q_1$ and $Q_2$ such that the pairs $(p,Q_1)$ and $(p,Q_2)$ lie within the admissible triangle \eqref{def: admissible}. Say
\begin{align*}
  \frac{1}{Q_i} = 
  \begin{cases} 
  \frac{1}{q_i}
  & \text{if } \frac{1}{p_i} \geq \frac{n-1}{2} \left( \frac{1}{2} - \frac{1}{q_i} \right), \\ 
  \frac{1}{2} - \frac{2}{n-1} \frac{1}{p_i}
  & \text{if } \frac{1}{p_i} \leq \frac{n-1}{2} \left( \frac{1}{2} - \frac{1}{q_i} \right).
  \end{cases}
\end{align*}
 Note that $1/Q_{i}\ge1/q_{i}$. By using Theorem \ref{thm: strichartz} with $\Theta_{i}=(n+1)(1/2-1/Q_i)/2$, we obtain
 \begin{align*}
  \|U\|_{L^{p_{1}}(\mathbb{R};\,H^{-\Theta_{1},Q_{1}}(G,\tau))}\,
  \lesssim\,
  \|U_0\|_{L^{2}(G,\tau)}\,+\,
  \|L\|_{L^{p_{2}'}(\mathbb{R};\,H^{\Theta_{2},Q_{2}'}(G,\tau))}.
  \end{align*}
Since $\theta_i-\Theta_i=n(1/Q_i-1/q_i)\ge0$, we deduce from Proposition \ref{prop: sobolev} that
\begin{align*}
  \|U\|_{L^{p_{1}}(\mathbb{R};\,H^{-\theta_{1},q_{1}}(G,\tau))}\,
  \lesssim\,
  \|U\|_{L^{p_{1}}(\mathbb{R};\,H^{-\Theta_{1},Q_{1}}(G,\tau))},
\end{align*}
and 
\begin{align*}
  \|L\|_{L^{p_{2}'}(\mathbb{R};\,H^{\Theta_{2},Q_{2}'}(G,\tau))}\,
  \lesssim\,
  \|L\|_{L^{p_{2}'}(\mathbb{R};\,H^{\theta_{2},q_{2}'}(G,\tau))}.
\end{align*}
Thus, the corollary is proved.
\end{proof}

\begin{remark}\label{rem: admissible 2}
  We present figures for the lower-dimensional cases $n=2$  and  $n=3$, which are of particular interest in studying the corresponding nonlinear problems. In both figures, the red region represents the admissible set \eqref{def: admissible} for the Strichartz inequality established via the $TT^*$ argument. Except for the red segment, the interior points are typically not available in the Euclidean setting. The blue region is obtained using the Sobolev embedding.

  \begin{figure}[h]
    \centering
    \begin{subfigure}{0.4\textwidth}
      \centering
      \begin{tikzpicture}
        \draw [->, line width=0.5] (2,0) -- (4,0);
        \draw [->, line width=0.5] (0,0) -- (0,4);

        \draw [blue, line width=0, opacity=0.2, fill=blue] (0,0) -- (0,2) -- (2,0) -- cycle;
        \draw [red, line width=0, opacity=0.2, fill=red] (0,2) -- (2,2) -- (2,0) -- cycle;

        \draw [red, line width=1] (2,0) -- (0,2);
        \draw [red, line width=1, dashed] (0,2) -- (2,2);
        \draw [red, line width=1, dashed] (2,2) -- (2,0);
        \draw [red,fill=red] (0,2) circle (0.05);
        \draw [red] (2,0) circle (0.05);

        \draw [blue, line width=1] (0,0) -- (0,2);
        \draw [blue, line width=1, dashed] (0,0) -- (2,0);
        \draw [blue] (0,0) circle (0.05);

        \node at (-0.3, -0.5) {$0$};
        \node at (4, -0.5) {$\frac{1}{p}$};
        \node at (-0.3, 3.8) {$\frac{1}{q}$};
        
        \node at (2, -0.5) {$\frac{1}{2}$};
        \node at (-0.3, 2) {$\frac{1}{2}$};
      \end{tikzpicture}
      \caption{Admissible set for $n=3$}
    \end{subfigure}
    \hfill
    \begin{subfigure}{0.4\textwidth}
      \centering
      \begin{tikzpicture}
        \draw [->, line width=0.5] (2,0) -- (4,0);
        \draw [->, line width=0.5] (0,0) -- (0,4);

        \draw [blue, line width=0, opacity=0.2, fill=blue] (0,0) -- (0,2) -- (1,0) -- cycle;
        \draw [red, line width=0, opacity=0.2, fill=red] (0,2) -- (2,2) -- (2,0) -- (1,0) -- cycle;

        \draw [red, line width=1] (1,0) -- (0,2);
        \draw [red, line width=1, dashed] (0,2) -- (2,2);
        \draw [red, line width=1, dashed] (2,2) -- (2,0);
        \draw [red, line width=1, dashed] (1,0) -- (2,0);
        \draw [red,fill=red] (0,2) circle (0.05);
        \draw [red] (2,0) circle (0.05);
        \draw [red] (1,0) circle (0.05);

        \draw [blue, line width=1] (0,0) -- (0,2);
        \draw [blue, line width=1, dashed] (0,0) -- (1,0);
        \draw [blue] (0,0) circle (0.05);

        \node at (-0.3, -0.5) {$0$};
        \node at (4, -0.5) {$\frac{1}{p}$};
        \node at (-0.3, 3.8) {$\frac{1}{q}$};
        
        \node at (2, -0.5) {$\frac{1}{2}$};
        \node at (-0.3, 2) {$\frac{1}{2}$};
      \end{tikzpicture}
      \caption{Admissible set for $n=2$}
    \end{subfigure}
  \end{figure}
\end{remark}

\noindent\textbf{Acknowledgement.}
The authors would like to thank Federico Cacciafesta for continuous discussion.
L.M. acknowledges support from the Deutsche
Forsch-ungsgemeinschaft (DFG, German Research Foundation) through SFB-TRR
352 -- Project 470903074 and from ERC CoG RAMBAS – Project-No. 10104-424. H.-W. Z receives funding from the Deutsche Forschungsgemeinschaft (DFG) through SFB-TRR 358/1 2023 ---Project 491392403. He would like to thank Zhipeng Yang for the warm hospitality during his visit to Yunnan Key Laboratory of Modern Analytical Mathematics and Applications (No. 202302AN360007).
J. Z. was supported by  National key R\&D program of China: 2022YFA1005700, National Natural Science Foundation of China(12171031) and Beijing Natural Science Foundation(1242011).


\printbibliography

\vspace{20pt}\address{
\noindent\textsc{Long Meng:}
\href{mailto: long.meng@lmu.de}{long.meng@lmu.de}\\
\textsc{Mathematisches Institut\\
  Ludwig-Maximilians-Universit\"at M\"unchen\\
 M\"unchen, Germany, 80333 }}

\vspace{10pt}\address{
\noindent\textsc{Hong-Wei Zhang:}
\href{mailto:zhongwei@math.upb.de}
{zhongwei@math.upb.de}\\
\textsc{Institut für Mathematik\\
Universität Paderborn\\
Paderborn, Germany, 33098 }

\vspace{10pt}\address{
\noindent\textsc{Junyong Zhang:}
\href{mailto: zhang\_junyong@bit.edu.cn}{zhang\_junyong@bit.edu.cn}\\
\textsc{Department of Mathematics\\
Beijing Institute of Technology\\
Beijing, China, 100081}

\end{document}